\documentclass[12pt]{article}
\usepackage{floatpag}

\usepackage{amsmath,amssymb,amstext}
\usepackage{graphicx}
\usepackage{geometry} 

\newcommand{\C}{\mathbb{C}}

\newcommand{\R}{\mathbb{R}}
\newcommand{\Z}{\mathbb{Z}}
\newcommand{\Digits}{\text{Digits}}

\title{Computing the moment polynomials of the zeta function}

\author{Michael O. Rubinstein
\footnote{Support for work on this paper was provided by the
National Science Foundation under awards DMS-0757627 (FRG grant),
and an NSERC Discovery Grant.} \phantom{a} and Shuntaro Yamagishi}

\begin{document}

\maketitle

\begin{abstract}

We describe a method to accelerate the numerical computation of the coefficients of
the polynomials $P_k(x)$ that appear in the conjectured asymptotics of the
$2k$-th moment of the Riemann zeta function. We carried out our method to compute
the moment polynomials for $k \leq 13$, and used these to experimentally
test conjectures for the moments up to height $10^8$.

\end{abstract}

\section{Introduction}

\label{section_intro}

For positive integer $k$, and any $\varepsilon > 0$,
Conrey, Farmer, Keating, Rubinstein, and Snaith conjectured~\cite{CFKRS} that
\begin{equation}
    \label{eq:moments zeta}
    \int_0^T |\zeta(1/2+it)|^{2k} dt= \int_0^T
    P_k\left(\log (t/(2\pi))\right) dt
    +O_{k,\varepsilon}\left( T^{1/2+\varepsilon} \right),
\end{equation}
with the constant in the $O$ term depending on $k$ and $\varepsilon$.

In the above equation $P_k$ is the polynomial of degree $k^2$ given implicitly by the
$2k$-fold residue
\newpage
\begin{eqnarray}
    \label{eq:P_k(x)}
     P_k(x)= && \frac{(-1)^k}{k!^2}\frac{1}{(2\pi i)^{2k}}
    \oint\cdots \oint \frac{G(z_1,
    \ldots,z_{2k})\Delta^2(z_1,\ldots,z_{2k})} {\displaystyle
    \prod_{l=1}^{2k} z_l^{2k}} \cr
    \cr
    \cr
    &&  \phantom{xxxxxxxxxxxxxxxxxxxxxxxxxxxx}\times
    e^{\tfrac x2\sum_{l=1}^{k}z_l-z_{l+k}}~dz_1\ldots dz_{2k} ,
\end{eqnarray}
with the path of integration over small circles about $z_l=0$, where
\begin{equation}
    \Delta(z_1,\ldots,z_{m})
    =\prod_{1\le i < j\le m}(z_j-z_i)
    = \left|  z_i^{j-1}  \right|_{m \times m}
\end{equation}
denotes the Vandermonde determinant,
\begin{equation}
    \label{eq:G}
    G(z_1,\ldots,z_{2k})= A_k(z_1,\ldots,z_{2k})
    \prod_{i=1}^k\prod_{j=1}^k\zeta(1+z_i-z_{j+k}) ,
\end{equation}
and $A_k$ is the Euler product
\begin{eqnarray}
    &&A_k(z_1,\ldots,z_{2k}) \notag \\
    &&=
    \label{eq:A_k}
    \prod_p
    \prod_{l,j=1}^k (1-p^{-1-z_l+z_{k+j}})
    \int_0^1 \prod_{j=1}^k
    \left(1-\frac{e(\theta)}{p^{\frac12 +z_j}}\right)^{-1}
    \left(1-\frac{e(-\theta)}{p^{\frac12 -z_{k+j}}}\right)^{-1}\,d\theta.
\end{eqnarray}
\newline Here $e(\theta) = \exp(2\pi i \theta)$.

We denote the coefficients of $P_k(x)$ by $c_r(k)$:
\begin{equation}
    \label{eq:c_r defn}
    P_k(x) =: \sum_{r=0}^{k^2} c_r(k) x^{k^2-r}.
\end{equation}

In order to arrive at this conjecture,
CFKRS considered a more general moment problem with
`shifts'. Because this general setting was central to our
computation, we describe their conjecture with shifts below.
Write the functional equation of zeta as
\begin{equation}
    \zeta(s) = \chi(s) \zeta(1-s),
\end{equation}
where
\begin{equation}
    \chi(s) := \pi^{s-1/2} \Gamma((1-s)/2)/\Gamma(s/2).
\end{equation}

It is a little more convenient to work with the Hardy $Z$-function, whose functional
equation and approximate functional equation are expressed more symmetrically than that of
the zeta function. It is defined as
\begin{equation}
    Z(s) = \chi(s)^{-1/2} \zeta(s),
\end{equation}
and satisfies
\begin{eqnarray}
    &&Z(s)=Z(1-s), \text{\,\, (because $\chi(s) \chi(1-s) = 1$)}, \notag \\
    &&Z(1/2+it) \in \R \text{ for $t \in \R$}, \notag \\
    &&|Z(1/2+it)| = |\zeta(1/2+it)|.
\end{eqnarray}

CFKRS took as their starting point the shifted moments:
\begin{equation}
    M(\alpha_1,\ldots,\alpha_{2k}) :=
    \int_{0}^{T} Z(1/2+it+\alpha_1)\cdots Z(1/2+it+\alpha_{2k})\,dt,
\end{equation}
where $\alpha_j \in \C$ are distinct and satisfy $-1/4 < \Re \alpha_j$.
When $\alpha=0$ the integrand is $|\zeta(1/2+it)|^{2k}$.

Substituting the approximate functional equation into each factor of the above integrand
\begin{equation}
    Z(s) = \chi(s)^{-1/2} \sum_{n\leq \sqrt{\frac{t}{2\pi}}} \frac{1}{n^s}
    + \chi(1-s)^{-1/2} \sum_{n\leq \sqrt{\frac{t}{2\pi}}} \frac{1}{n^{1-s}}
    +O(t^{-\sigma/2}),
\end{equation}
$s=\sigma+it$, $0 < \sigma < 1$,
CFKRS applied the following {\it heuristic} steps:
\begin{itemize}
    \item[a)] Ignore the $O(t^{-\sigma/2})$ and expand the product to get $2^{2k}$ terms, each a product
              $2k$ sums.
    \item[b)]
    Of the $2^{2k}$ terms, only the terms with the same number of $s$'s and $1-s$'s
    contribute to the asymptotics. Reasoning: $\chi(s)$ is highly oscillatory, so
    cancellation occurs unless each $s$ gets paired with a $1-s$.
    \item[c)]
    For any such term, only the diagonal (`$m_1 m_2 \ldots m_k = n_1 n_2 \ldots n_k$')
    contributes when the sums are multiplied out.
    \item[d)]
    Extend the truncated diagonal sums to infinity, replacing the sums that diverge
    with their analytic continuation (the assumption we stated
    earlier, $-1/4 < \Re \alpha_j$, is used when obtaining
    the analytic continuation of the diagonal sums).
\end{itemize}
The steps in this heuristic recipe are not justifiable, and individually
are not even true! The terms that are dropped cannot be neglected, and it appears that
some sort of cancellation takes place amongst these terms so that,
in the end, the above steps do apparently result in a correct conjecture,
as described below.

Let
\begin{equation}
    \label{eq:H a}
    H(z_1,\ldots,z_{2k};x) :=
    \exp\left( \frac{x}{2} \sum_1^k z_j-z_{j+k}\right)
    G(z_1,\ldots,z_{2k}).
\end{equation}
The first conjectured asymptotic formula of CFKRS, which we refer to as the combinatorial sum,
for shifted moments reads:
\begin{equation}
    M(\alpha_1,\ldots,\alpha_{2k}) \sim
    \int_{0}^{T} P_k\left(\alpha,\log \tfrac{t}{2 \pi}\right) ~dt,
\end{equation}
where
\begin{equation}
    P_k(\alpha,x) =
    \sum_{\sigma \in \Xi} H( \alpha_{\sigma(1)},\ldots,\alpha_{\sigma(2k)};x),
    \label{eq:comb sum}
\end{equation}
and $\Xi \subset S_{2k}$ is the set of $\binom{2k}{k}$ permutations
such that $\sigma(1)<\cdots < \sigma(k)$ and $\sigma(k+1)<\cdots < \sigma(2k)$.
The terms in this set correspond to the number of ways to select, from approximate
functional equation for $Z(s)$, the same number $k$ of $s$'s and $(1-s)$'s.

CFKRS, also expressed their sum of $\binom{2k}{k}$ terms
as a $2k$-fold residue. We reproduce their second formula for the shifted moments:
\begin{eqnarray}
    \label{eq:P_k alpha}
    P_k(\alpha,x)=
    \frac{(-1)^k}{k!^2}\frac{1}{(2\pi i)^{2k}}
    \oint \cdots \oint
    &&\frac{G(z_1,
    \ldots,z_{2k})\Delta(z_1,\ldots,z_{2k})^2} {\displaystyle
    \prod_{l=1}^{2k} \prod_{j=1}^{2k}(z_l-\alpha_j)}
    \\
    &&\times e^{\tfrac
    x2\sum_{l=1}^{k}z_l-z_{l+k}}dz_1\ldots dz_{2k}. \notag
\end{eqnarray}
The Vandermonde in the above residue vanishes whenever any of the $z_j$'s are
equal, thus restricting the residues to contributions from choosing a distinct
$\alpha$ for each $z$, i.e. $z_j=\alpha_{\sigma(j)}$, where $\sigma$ is a
permutation of $1,2,\ldots,2k$. The symmetry of the function $G$ with respect
to the first $k$ variables, and with respect to the last $k$ variables allows
one to write the resulting residues as a sum over the $\binom{2k}{k}$
permutations $\Xi$ in~\eqref{eq:comb sum}. The residue computation establishing
the equivalence of~\eqref{eq:P_k alpha} and~\eqref{eq:comb sum} is carried out
in Lemma 2.5.3 of~\cite{CFKRS}.

While the $\binom{2k}{k}$ terms of the combinatorial sum~\eqref{eq:comb sum} for $P_k(\alpha,x)$
each have poles of order $k^2$ at $\alpha=0$, the above is analytic in a neighbourhood of $\alpha=0$
which shows that these poles must cancel. Working with shifts allowed CFKRS to get beyond these poles.

To get formula \eqref{eq:P_k(x)} for $P_k(x)$, set $\alpha=0$ in~\eqref{eq:P_k
alpha}. Even though the formula conjectured by CFKRS is complicated, it does
seem to correctly predict the moments. Some evidence, both theoretical and
numerical, in its favour was presented in~\cite{CFKRS} and~\cite{CFKRS2}. For
instance, the formulas predicted by CFKRS match with known theorems, including
lower terms, for $k=1,2$. For $k=1$, the full asymptotics were obtained by
Ingham~\cite{I}, when $|\alpha_1|,|\alpha_2|<1/2$. For $k=2$, the full
asymptotics was proven by Heath-Brown~\cite{H-B} when $\alpha_j=0$, and
Motohashi with shifts $\alpha_j$ in some neighbourhood of 0.

Perhaps the best evidence for these heuristics is numerical. CFKRS computed the
moment polynomials for $k\leq 7$, and tested the moment conjecture for $T$
roughly of size $10^6$. In this paper we test the conjecture for $k\leq 13$ and
$T$ up to $10^8$. Additional numerics are discussed in~\cite{HO}.

The heuristics do not shed much light on the nature of the remainder term. The
related multiple Dirichlet series approach of Diaconu, Goldfeld, and
Hoffstein~\cite{DGH} does predict, based on conjectured analytic properties of
the relevant multiple Dirichlet series, a remainder term of size
$O_{k,\varepsilon}\left( T^{1/2+\varepsilon} \right)$.

In~\cite{CFKRS2}, two methods are described for computing the coefficients $c_r(k)$
of $P_k(x)$. Their first method, based on deriving explicit formulas
for the $c_r(k)$, taking  equation~\eqref{eq:P_k(x)} as the starting point,
is useful for obtaining the coefficients to very
high precision, but is limited to smaller values of $r$. In~\cite{CFKRS}\cite{CFKRS2},
their first method was used to compute $c_r(k)$ for $r\leq 9$ and various values of $k$.

Their second method is more useful for computing all $k^2+1$ coefficients of
$P_k(x)$, and is easier to implement, but is more limited in the amount of precision
that it can achieve. It uses~\eqref{eq:comb sum}, small shifts $\alpha_j$, and high
precision to capture cancellation resulting from high order poles that cancel in the
terms of this sum.

In this paper we describe a `cubic accelerant' variant of the second method
presented in~\cite{CFKRS2} for computing coefficients of the moment
polynomials. This allowed us to extend tables of the $c_r(k)$ for $k\leq 13$,
and, $0 \leq r \leq k^2$. We also computed many of the coefficients to greater
precision. Finally, we used our tables of coefficients to test the moment
conjectures up to $k\leq 13$ and $T$ up to $10^8$.

Going up to $k=13$, using~\eqref{eq:comb sum}, or more
precisely using~\eqref{eq:c_r as limit}, is substantial, because computing the
coefficients to $D$ digits accuracy involves evaluating $k^2$ sums (one for
each $r$), each sum involving $2k \choose k$ terms ($10,400,600$ for $k=13$),
with working precision of roughly $D\times k^2$ digits accuracy. For example,
about $2000$ digits are required for $k=13$ with desired precision of $12$
digits. The process is made even more challenging, by the fact that each term
involves a complicated infinite multivariate product over primes. The computer
data storage (RAM) requirements are also large because of the number of terms, $2k
\choose k$, that we are computing/updating one prime at a time to high
precision.

The basic idea of the accelerant is to approximate the tail of the Euler
product~\eqref{eq:A_k} defining $A_k$, with an expression, given
in~\eqref{eq:cubic acc} that accounts for all terms of degree $\leq 3$ in
$1/p$. In our application, we will take the $z_j$'s appearing in this
approximation to be very small. This approximation gives, for large $p$, the
local factor of $A_k$ up to a remainder of size roughly $O(k^8 p^{-4})$, so
that the overall contribution, from all $p$ greater than a given large $P$, is
approximated up to a remainder roughly of size $O(k^8 P^{-3} \log(P)^{-1})$.
This is described in Section~\ref{sec:choice P}.

For instance, the use of our cubic approximation~\eqref{eq:cubic acc} allowed
us to truncate the multivariate Euler product $A_{13}$ at $P=1699$, and achieve
the quality of coefficients described in Table~\ref{table:coeff_13} of
Section~\ref{sec:implementation}. Without this approximation, we would have
needed to truncate at $P \approx 10^9$. For $k=4$, the results we obtained with
$P=942,939,827$ would have required $P \approx 10^{27}$ without any
acceleration.

\section{Computing the coefficients of $P_k(x)$}
\label{sec:method 1}


The first method of CFKRS for computing the coefficients of $P_k(x)$
involved expanding, in a multivariate Taylor series, the integrand
of~\eqref{eq:P_k(x)}, and working out a technique for expressing the resulting
residue (for a general $k$), giving formulas for $c_r(k)$.
For example,
\begin{eqnarray}
    \label{eq:c0 c1 c2}
    c_0(k) &=& a_k \prod_{j=0}^{k-1}\frac{j!}{(j+k)!} \notag \\
    c_1(k) &=& c_0(k) 2 k^2 (\gamma_0 k + B_k(1;)) \notag \\
    c_2(k) &=& c_0(k) k^2 (k-1)(k+1) \\
           &&\times (
               2(B_k(1;)+\gamma_0 k)^2 - \gamma_0^2 -2\gamma_1+B_k(1,1;)-B_k(1;1)
           ),
           \notag
\end{eqnarray}
where
\begin{equation}
\label{eq:a_k}
    a_k= \prod_p \left(1-p^{-1}\right)^{k^2} {}_2F_1(k,k;1;1/p),
\end{equation}
the $\gamma_j$'s defined by
\begin{equation}
    \label{eq:zeta series}
    s \ \zeta(1+s) = 1 + \gamma_0 s - \gamma_1 s^2 +
    \frac{\gamma_2}{2!} s^3 -\frac{\gamma_3}{3!} s^4 + \cdots,
\end{equation}
and
\begin{eqnarray}
    \label{eq:B_k}
    B_k(1;) &=& \sum_p
    \frac{k\log(p)}{p-1} -
    \frac{\log(p) k\,{}_2F_1(k+1,k+1;2;1/p)}{p\,{}_2F_1(k,k;1;1/p)}
    \\
    B_k(1,1;) &=& -\sum_p
    \left(
         \frac{\log(p)^2 k^2{}_2F_1(k+1,k+1;2;1/p)^2}{p^2\,{}_2F_1(k,k;1;1/p)^2} -
         \frac{\log(p)^2{k+1\choose 2}{}_2F_1(k+2,k+2;3;1/p)}{p^2\,{}_2F_1(k,k;1;1/p)}
    \right)
    \notag \\
    B_k(1;1) &=& \sum_p
    \frac{p\log(p)^2}{(p-1)^2} +
    \biggl(
         \frac{\log(p)^2 k^2{}_2F_1(k+1,k+1;2;1/p)^2}{p^2 {}_2F_1(k,k;1;1/p)^2} - \notag \\
         &&\frac{\log(p)^2{}_2F_1(k+1,k+1;1;1/p)}{p\,{}_2F_1(k,k;1;1/p)}
    \biggr)
    \notag \\
    B_k(2;) &=& -\sum_p
    \frac{kp\log(p)^2}{(p-1)^2} +
    \biggl(
         \frac{\log(p)^2 k^2{}_2F_1(k+1,k+1;2;1/p)^2}{p^2\,{}_2F_1(k,k;1;1/p)^2} - \notag \\
         &&\frac{\log(p)^2{k+1\choose 2}\,{}_2F_1(k+2,k+2;3;1/p)}{p^2\,{}_2F_1(k,k;1;1/p)} -
         \frac{\log(p)^2 k\,{}_2F_1(k+2,k+1;2;1/p)}{p\,{}_2F_1(k,k;1;1/p)}
    \biggr),
    \notag
\end{eqnarray}
with ${_2F_1}$ Gauss' hypergeometric function.

These formulas quickly get much more complicated. In practice CFKRS were able to
use this method for $r \leq 9$ and compute numerical approximations for
all the coefficients of, for example, $P_3(x)$. One advantage of these formulas,
expressed as sums over primes, is that one can apply Mobius inversion to accelerate
the convergence of these sums, and obtain high precision values of the coefficients.

It soon became apparent~\cite{CFKRS} from numerical values of $c_r(k)$ that the
leadings coefficients of $P_k(x)$, i.e. associated to the larger powers of $x$,
are very small in comparison to the lower terms. Thus, in order to
meaningfully test the moment conjecture for zeta, which involves the moment
polynomial evaluated at the slowly growing function $\log(t/2\pi)$ (and this
hardly changes over the range of $t$ in which we can gather significant data
for $\zeta(1/2+it)$), one needs many coefficients of the moment polynomials.
See also~\cite{HR} which discusses the uniform asymptotics of these coefficients.

Consequently, a second practical method, relying on the combinatorial
sum~\eqref{eq:comb sum}, was developed for computing numerical approximations for {\it all}
$k^2$ coefficients of the moment polynomial $P_k(x)$.

We detail our computational approach, implementation, and numerical results in
the next two sections.

\section{Our numerical evaluation of $c_r(k)$}

\label{section_numerical}

The polynomial $P_k(x)$ given by~(\ref{eq:P_k(x)}) is the special case
$\alpha_1 = \ldots = \alpha_{2k}= 0$ of the function $P_k(\alpha,x)$
in~(\ref{eq:P_k alpha}). CFRKS's second method for computing the coefficients
of $P_k(x)$ relies on their equation~(\ref{eq:comb sum}) for $P_k(\alpha,x)$.
However, the terms in~(\ref{eq:comb sum}) have poles if the $\alpha_j$'s are
not distinct, coming from the product of $k^2$ zetas,
\begin{equation}
    \label{eq:zetas}
    \prod_{i=1}^k\prod_{j=1}^k\zeta(1+z_i-z_{j+k}),
\end{equation}
that appear in the function $G$. So we cannot simply substitute $\alpha_j = 0$.

Instead we take the limit as $\alpha_j \to 0$ while making sure that all the
$\alpha_j$'s are distinct. Because of the poles, each individual term
in~\eqref{eq:comb sum}  becomes very large when $\alpha$ is small, and high
precision is needed to see one's way through the resulting cancellation of the
poles as we sum across the $2k \choose k$ terms of the combinatorial sum.

More precisely, consider
\begin{equation}
    H(z_1,\ldots,z_{2k};x) =
    \exp\left( \frac{x}{2} \sum_1^k z_j-z_{j+k}\right)
    A_k(z_1,\ldots,z_{2k})
    \prod_{i=1}^k\prod_{j=1}^k\zeta(1+z_i-z_{j+k}),
\end{equation}
and let
\begin{equation}
    \delta_j = j \delta,
\end{equation}
where $\delta \in {\mathbb C}$ is a small number. In practice $\delta$ was
of the form $10^{-D}$ for some positive integer $D$.

Using (\ref{eq:comb sum}) we obtain
\begin{equation}
    \label{eq:P_k as limit}
    P_k(x) =
    \lim_{\delta \to 0}
    \sum_{\sigma \in \Xi} H( \delta_{\sigma(1)},\ldots,\delta_{\sigma(2k)};x).
\end{equation}

As in~\cite{CFKRS2}, we expand $\exp\left( \frac{x}{2} \sum_1^k
z_j-z_{j+k}\right)$ in its Taylor series, and pull out the coefficient of
$x^{k^2-r}$, to get
\begin{equation}
    \label{eq:c_r as limit}
    c_r(k) =
    \frac{1}{2^{k^2-r} (k^2-r)!}
    \lim_{\delta \to 0}
    \sum_{\sigma \in \Xi} H_r( \delta_{\sigma(1)},\ldots,\delta_{\sigma(2k)}),
\end{equation}
where
\begin{equation}
    \label{eq:H_r}
    H_r(z_1,\ldots,z_{2k}) =
    \left(\sum_1^k z_j-z_{j+k}\right)^{k^2-r}
    A_k(z_1,\ldots,z_{2k})
    \prod_{i=1}^k\prod_{j=1}^k\zeta(1+z_i-z_{j+k}).
\end{equation}
Notice that, as a function of $\delta$,  $H_r( \delta_{\sigma(1)},\ldots,\delta_{\sigma(2k)})$
has a pole at $\delta=0$ of order $r$, because the first factor above cancels $k^2-r$ of
the $k^2$ poles of the double product of zetas. These poles must cancel when summed
over permutations $\sigma$, otherwise we would not obtain the lhs,
$c_r(k)$, as $\delta \to 0$. Therefore, because the sum in~\eqref{eq:c_r as limit}
is analytic about $\delta=0$, we may write
\begin{equation}
    \label{eq:c_r as limit b}
    c_r(k) =
    \frac{1}{2^{k^2-r} (k^2-r)!}
    \left(
        \left(
            \sum_{\sigma \in \Xi} H_r(
            \delta_{\sigma(1)},\ldots,\delta_{\sigma(2k)})
        \right)
        +O(|\delta|)
    \right),
\end{equation}
with the implied constant in the remainder term depending on $k$ and $r$.
In our implementation, we neglected the contribution from the $O(|\delta|)$ term for
reasons that are described at the end of~\ref{sec:choice P}.

One complication in evaluating the above for a given $k$ and $\delta$
is that $A_k(z_1,\ldots,z_{2k})$ is expressed as an infinite product over
primes as described by~(\ref{eq:A_k}).

While CFKRS used a `quadratic accelerant' for evaluating the multivariate Euler
product, we implemented a cubic accelerant. This has the advantage of allowing
us to truncate the Euler product sooner.

To evaluate the convergent product $A_k(z_1,\ldots,z_{2k})$, we break up the
product over primes into $p \leq P$ and $p>P$, where $P$ is a large number. For
the first portion $p\leq P$, we use the following identity, derived in Section
2.6 of~\cite{CFKRS},
\begin{equation}
    \label{eq:A_k euler product 2}
    A_k(z_1,\ldots,z_{2k}) =
    \prod_p \sum_{j=1}^k
    \prod_{i\neq j}
        \frac{\displaystyle \prod_{m=1}^k (1-p^{-1+z_{i+k}-z_m})}
        {1-p^{z_{i+k}-z_{j+k}}},
\end{equation}
to numerically compute the local factor of $A_k(z_1,\ldots,z_{2k})$ for specific values of
$p$ and $z_1,\ldots,z_{2k}$.
Ideally, we would compute the local factor for all $p$ using this formula, but,
because there are infinitely many primes, we must eventually stop. However, we
will describe a method to approximate the contribution from all $p>P$, thus
allowing us to attain higher precision in our computation with fewer primes.
The choice of $P$ is described in Section~\ref{sec:choice P}.

Some care must be taken to account for the fact that individual terms
in~\eqref{eq:A_k euler product 2} also have poles. While these poles cancel out
when summed over $j$, see the paragraph following equation (2.6.16)
in~\cite{CFKRS}, they cause some additional loss of precision in our
application. We are evaluating $A_k(z_1,\ldots,z_{2k})$ at distinct, but small
values of $z_j$. Therefore, when evaluating the sum over $j$, additional
cancellation and hence loss of precision occurs affecting the leading $(k-1) D$
digits of the truncated Euler product for $A_k(z_1,\ldots,z_{2k})$, where $z_j
\approx 10^{-D}$, from the poles of order $k-1$ of the individual terms summed
in~\eqref{eq:A_k euler product 2}.

For the contribution of the second portion $p>P$, we approximate each local
factor appearing in~(\ref{eq:A_k}) by a product of zeta functions that captures
the terms up to degree three in $1/p$ of its multivariate Dirichlet series.
A cubic approximation can be obtained by first substituting
$u_j=p^{-1/2 - z_j}$ and $w_j=p^{-1/2 + z_{k+j}}$ into the local factor
of~(\ref{eq:A_k}),
\begin{equation}
    \label{eq:uw integral}
    \prod_{l,j=1}^k (1-u_lw_j)
    \int_0^1 \prod_{j=1}^k
    (1-u_je(\theta))^{-1}
    (1-w_je(-\theta))^{-1}\,d\theta,
\end{equation}
and then working out the terms, in the multivariate Maclaurin series,
up to degree six, in $u_j$ and $w_j$.

Notice that the integral over $\theta$ pulls out just the terms with the same
number of $u$'s and $w$'s. This results in monomials only of even
degree appearing. The integral of any other term, which does not have the same
number of $u$'s and $w$'s, is zero because it contains a non-zero integer power
of $e(\theta)$.

Also, observe that the local factor of~(\ref{eq:uw integral}) is symmetric in the $u$'s
and, separately in the $w$'s, meaning if the $u_l$'s are permuted
the expression remains invariant and similarly for the $w_l$'s. Also, it is
symmetric with $u$ and $w$, i.e. if all the $u$'s and $w$'s are swapped the
expression remains the same.

Therefore, to get terms up to degree six, we can determine the coefficients of
representative terms involving $u_1,u_2,u_3$ and
$w_1,w_2,w_3$, and then symmetrize the resulting expressions over all the $u$'s
and $w$'s. More precisely, to get all terms of degrees 2, 4, and 6, it is
sufficient to consider only the monomials: $u_1 w_1$, $u_1 u_2 w_1 w_2$,
$u_1 u_2 w_1^2$, $u_1^2 w_1^2$, $u_1u_2u_3w_1w_2w_3$, $u_1u_2u_3w_1^{2}w_2$,
$u_1u_2u_3w_1^{3}$, $u_1^{2}u_2w_1^{2}w_2$, $u_1^{2}u_2w_1^{3}$, and
$u_1^{3}w_1^{3}$ instead of every possible monomial, and then exploit
symmetry.

Finally, we can simplify further. The integral over $\theta$ simply plays the
role of pulling out terms with the same number of $u$'s and $w$'s. So, instead
of~\eqref{eq:uw integral}, we can work more directly with the function
\begin{equation}
    \label{eq:uw function}
    \prod_{l,j=1}^k (1-u_lw_j)
    \prod_{j=1}^k
    (1-u_j)^{-1}
    (1-w_j)^{-1}.
\end{equation}
The multivariate Maclaurin series of the above coincides with that of~\eqref{eq:uw integral}
for those terms that have the same number of $u$'s and $w$'s. Furthermore,
because we are focusing just on terms involving $u_1,u_2,u_3$ and $w_1,w_2,w_3$,
we can set $u_j=w_j=0$ for all $j\geq 4$. Finally, to get terms up to degree six
with the same number of $u$'s and $w$'s, we can expand each factor in the
denominator as a geometric series of degree 3. We therefore consider:
\begin{equation}
    \label{eq:uw function b}
    \prod_{l,j=1}^3 (1-u_lw_j)
    \prod_{j=1}^3
    (1+u_j+u_j^2+u_j^3)
    (1+w_j+w_j^2+w_j^3).
\end{equation}
We expand out the above and
tabulate, in Table~\ref{table:uw}, the coefficients for representative monomials
with the same number of $u$'s and $w$'s, up to degree 6, in the multivariate
Maclaurin series of the above function, and, hence, equivalently, in~\eqref{eq:uw integral}.

\begin{table}[h]
\centerline{
\begin{tabular}{|c|c|}
\hline
monomial  & coefficient \\
\hline
$u_1 w_1$ & 0 \\
$u_1 u_2 w_1 w_2$ & -1 \\
$u_1 u_2 w_1^2$ & 0 \\
$u_1^2 w_1^2$ & 0 \\
$u_1u_2u_3w_1w_2w_3$ & 4 \\
$u_1u_2u_3w_1^{2}w_2$ & 1 \\
$u_1u_2u_3w_1^{3}$ & 0 \\
$u_1^{2}u_2w_1^{2}w_2$ & 0 \\
$u_1^{2}u_2w_1^{3}$ & 0 \\
$u_1^{3}w_1^{3}$ & 0 \\
\hline
\end{tabular}
}
\caption{The second column lists the coefficients that
appear with representative monomials, up to degree 6, in the multivariate
Maclaurin expansion of~\eqref{eq:uw integral}.}
\label{table:uw}
\end{table}

Therefore, symmetrizing, we have that~\eqref{eq:uw integral} equals
\begin{eqnarray*}
1-\sum_{{1 \leq i_1<i_2 \leq k}\atop {1 \leq j_1<j_2 \leq k}}
    u_{i_1}u_{i_2}w_{j_1}w_{j_2} + \phantom{2} 4\sum_{{1 \leq i_1<i_2<i_3 \leq k}\atop {1 \leq j_1<j_2<j_3 \leq k}}
    u_{i_1}u_{i_2}u_{i_3}w_{j_1}w_{j_2}w_{j_3}
\end{eqnarray*}
\begin{eqnarray}
    + \sum_{{1 \leq i_1<i_2<i_3 \leq k}\atop {1 \leq j_1 \neq j_2\leq k}}
    u_{i_1}u_{i_2}u_{i_3}w_{j_1}^{2}w_{j_2}+ \sum_{{1 \leq i_1 \neq i_2\leq k}\atop {1 \leq j_1<j_2<j_3 \leq k}}
    u_{i_1}^{2}u_{i_2}w_{j_1}w_{j_2}w_{j_3} + \ldots
\end{eqnarray}

Undoing the substitution for $u_l$'s and $w_j$'s, gives the following expansion for the
local factors in~\eqref{eq:A_k}:
\begin{eqnarray}
\label{eq:local factor deg 3}
    1 &&-  \sum_{{1 \leq i_1<i_2 \leq k}\atop {1 \leq j_1<j_2 \leq k}} p^{-2-z_{i_1}-z_{i_2}+z_{k+j_1}+z_{k+j_2}}\cr
    \cr
    &&+ \phantom{x} 4\sum_{{1 \leq i_1<i_2<i_3 \leq k}\atop {1 \leq j_1<j_2<j_3 \leq k}} p^{-3-z_{i_1}-z_{i_2}-z_{i_3}+z_{k+j_1}+z_{k+j_2}+z_{k+j_3}}\cr
    \cr
    &&\phantom{1234567890}+ \sum_{{1 \leq i_1<i_2<i_3 \leq k}\atop {1 \leq j_1 \neq j_2\leq k}} p^{-3-z_{i_1}-z_{i_2}-z_{i_3}+2z_{k+j_1}+z_{k+j_2}}\cr
    \cr
    && \phantom{123456789012345678} + \sum_{{1 \leq i_1 \neq i_2\leq k}\atop {1 \leq j_1<j_2<j_3 \leq k}} p^{-3-2z_{i_1}-z_{i_2}+z_{k+j_1}+z_{k+j_2}+z_{k+j_3}}
    \notag \\
    && \phantom{12345678901234567812345678}+\ldots,
\end{eqnarray}
which we then approximate by the following product:
\begin{eqnarray}
    \label{eq:prod}
    &&\prod_{{1 \leq i_1<i_2 \leq k}\atop {1 \leq j_1<j_2 \leq k}}
    (1-p^{-2-z_{i_1}-z_{i_2}+z_{k+j_1}+z_{k+j_2}}) \cr
    && \phantom{12345} \times \prod_{{1 \leq i_1<i_2<i_3 \leq k}\atop {1 \leq j_1<j_2<j_3 \leq k}}
    (1-p^{-3-z_{i_1}-z_{i_2}-z_{i_3}+z_{k+j_1}+z_{k+j_2}+z_{k+j_3}})^{-4} \cr
    && \phantom{12345678901} \times \prod_{{1 \leq i_1<i_2<i_3 \leq k}\atop {1 \leq j_1 \neq j_2\leq k}}
    (1-p^{-3-z_{i_1}-z_{i_2}-z_{i_3}+2z_{k+j_1}+z_{k+j_2}})^{-1} \cr
    && \phantom{12345678901234567} \times \prod_{{1 \leq i_1 \neq i_2\leq k}\atop {1 \leq j_1<j_2<j_3 \leq k}}
    (1-p^{-3-2z_{i_1}-z_{i_2}+z_{k+j_1}+z_{k+j_2}+z_{k+j_3}})^{-1}.
    \notag \\
\end{eqnarray}
The last step can be seen by expanding each factor in a geometric series and
comparing the terms, up to those containing a $1/p^3$, with those
in~\eqref{eq:local factor deg 3}.
We also remark that, had we wanted a quartic
approximation, then slightly more care would be needed
as the first product above would, on expanding in geometric series,
interact with the quartic terms.

The product in~\eqref{eq:prod} allows us to approximate
the tail, i.e. for $p>P$, of $A_k(z_1,\ldots,z_k)$ in
terms of the Riemann zeta function:
$$
    \prod_{p>P}
    \prod_{l,j=1}^k (1-p^{-1-z_l+z_{k+j}})
    \int_0^1 \prod_{j=1}^k
    \left(1-\frac{e(\theta)}{p^{\frac12 +z_j}}\right)^{-1}
    \left(1-\frac{e(-\theta)}{p^{\frac12 -z_{k+j}}}\right)^{-1}\,d\theta
$$
$$
    \approx \frac{
        \prod_{{1 \leq i_1<i_2 \leq k}\atop {1 \leq j_1<j_2 \leq k}}\zeta(2+z_{i_1}+z_{i_2}-z_{k+j_1}-z_{k+j_2})^{-1}
    }
    {
         \prod_{p\leq P} \prod_{{1 \leq i_1<i_2 \leq k}\atop {1 \leq j_1<j_2 \leq k}} (1-p^{-2-z_{i_1}-z_{i_2}+z_{k+j_1}+z_{k+j_2}})
    }\phantom{12345678901234567890123456789012345678901234567890}
$$
$$
    \times \frac{
         \prod_{{1 \leq i_1<i_2<i_3 \leq k}\atop {1 \leq j_1<j_2<j_3 \leq k}}\zeta(3+z_{i_1}+z_{i_2}+z_{i_3}-z_{k+j_1}-z_{k+j_2}-z_{k+j_3})^{4}
    }
    {
        \prod_{p\leq P} \prod_{{1 \leq i_1<i_2<i_3 \leq k}\atop {1 \leq j_1<j_2<j_3 \leq k}} (1-p^{-3-z_{i_1}-z_{i_2}-z_{i_3}+z_{k+j_1}+z_{k+j_2}+z_{k+j_3}})^{-4}
    } \phantom{12345}
$$
\begin{equation}
    \label{eq:cubic acc}
    \phantom{1234} \times
    \frac{
         \prod_{{1 \leq i_1<i_2<i_3 \leq k}\atop {1 \leq j_1 \neq j_2\leq k}}\zeta(3+z_{i_1}+z_{i_2}+z_{i_3}-2z_{k+j_1}-z_{k+j_2})
    }
    {
         \prod_{p\leq P} \prod_{{1 \leq i_1<i_2<i_3 \leq k}\atop {1 \leq j_1 \neq j_2\leq k}} (1-p^{-3-z_{i_1}-z_{i_2}-z_{i_3}+2z_{k+j_1}+z_{k+j_2}})^{-1}
    }
\end{equation}
$$
    \phantom{1234567890123456789012} \times \frac{
         \prod_{{1 \leq i_1 \neq i_2\leq k}\atop {1 \leq j_1<j_2<j_3 \leq k}}\zeta(3+2z_{i_1}+z_{i_2}-z_{k+j_1}-z_{k+j_2}-z_{k+j_3})
    }
    {
         \prod_{p\leq P} \prod_{{1 \leq i_1 \neq i_2\leq k}\atop {1 \leq j_1<j_2<j_3 \leq k}} (1-p^{-3-2z_{i_1}-z_{i_2}+z_{k+j_1}+z_{k+j_2}+z_{k+j_3}})^{-1}
    }.
$$
\section{Implementation and tables of coefficients $c_r(k)$}
\label{sec:implementation}

Our code was implemented  in {\tt C++} using the GNU MPFR library~\cite{FHLPZ},
along with Jon Wilkening's {\tt C++} wrapper for MPFR. MPFR is based on
GMP, the GNU multiprecision library. We used gcc, the GNU C compiler, to compile
our code, with the `-fopenmp' option in order to enable the use of
OpenMP directives in our code. This allowed us to carry out some of the
key steps in parallel for a given $k$, using several cores of our machine.
Computations were carried out on an SGI Altix 3700 computer
with 64 Itanium2 processors and 192 GB of shared memory.

For each $k$, we selected a precision, specified by the number of
digits desired, `Digits', for the final output, and let
$\delta=10^{-\Digits}$. For example, we used Digits$=25$, i.e.
$\delta=10^{-25}$, for $k=4$.
We then put $\delta_j = j \delta$, $1 \leq j \leq k$,
and set about computing the sum~\eqref{eq:c_r as limit b}, using
our cubic multivariate approximation, i.e. equation~\eqref{eq:cubic acc},
for the tail of the infinite product $A_k(z_1,\ldots,z_{2k})$.

Observe, in~\eqref{eq:H_r}, that the dependence on $r$ manifests only at
the factor:
\begin{equation}
    \left(\sum_1^k z_j-z_{j+k}\right)^{k^2-r}.
\end{equation}
Therefore, we were able to store and recycle all the other quantities across
$0 \leq r \leq k^2$.

We record one important hack that we used several times in our program. While the
double product of zetas in~\eqref{eq:H_r} involves $k^2$ factors, many of these are
repeated since there are just $4k-2$ possible values of $\delta_a
-\delta_b=(a-b) \delta$, where $a,b$ are distinct integers in $[1,2k]$. The
same holds, for each $p$, in the double products in~\eqref{eq:A_k euler product 2}.

Likewise, while the products in~\eqref{eq:cubic acc} involve up to $O(k^6)$ factors,
these appear with multiplicity, and there are just $O(k)$ distinct factors. This is
true both for the product of zetas in the numerator, and also, for each $p$, the factors
that appear in the denominator.

We exploited these multiplicities by computing and storing a table of the
distinct values of zeta that appear and, for each $p$, of the distinct factors that
occur. Furthermore, we took advantage of the fact that the powers of $p$ that
occur, other than $1/p$, are of the form $p^{m \delta}$, where $m \in \Z$,
and thus computed and stored them by repeated multiplication of $p^{\delta}$
and of $1/p^{\delta}$.

This allowed us to avoid recomputing the same quantities
repeatedly and also to cut back significantly on the amount of high precision
multiplications needed.
For instance, we looped through the various indices in~\eqref{eq:cubic acc} to
count which factors appear with which multiplicities. This could be done
simply, for our range of $k$, using 32 bit integer arithmetic, with $O(k)$
exponentiations and multiplications then carried out for each $p$
in~\eqref{eq:cubic acc}, rather than $O(k^6)$ multiplications.

To account for the high amount of cancellation that occurs
as a consequence of the poles of the individual terms in~\eqref{eq:c_r as limit b},
we let our working precision be equal to
\begin{equation}
    \label{eq:working precision}
    \text{WorkingDigits} = (k^2+k-1+6) \Digits,
\end{equation}
and carried out our computations using these many digits.
The $k^2$ was to account for the largest order poles occurring in $H_r$,
when $r=k^2$, of order $k^2$. While we could have gotten away with less
precision for smaller $r$, we recycled most of the computed quantities
across all $r$.
The $k-1$ accounts for cancellation amongst the poles of the terms in the
sum over $j$ in~\eqref{eq:A_k euler product 2}. Finally we needed to have some working
precision left over, after all the cancellation, to capture $c_r(k)$ to
Digits precision. The $+6 \Digits$ was chosen to give us some leeway.
For example, we had WorkingDigits$=625$ for $k=4$ and Digits$=25$, and
WorkingDigits$=2244$ for $k=13$ and Digits$=12$.

While it would have been preferable to use a larger value of Digits for all our
$k$, rather than, for instance, a smaller value of Digits for $k=13$, two things made
this prohibitive.

The first was computing time. For $k=13$, our program ran on 6 processors
(using OpenMP to parallelize our code) for around 6 months, representing about
3 CPU years. Setting Digits$=25$ for $k=13$, say, and thus WorkingDigits of
around 4000 rather than 2000, would have at least doubled the amount of
computing time needed to carry out each arithmetic operation, as explained in
the next two subsections.
Furthermore, to
see the benefit of using a smaller value of $\delta=10^{-25}$ for $k=13$ would
have required us to truncate our product over primes at a much larger $P$,
around $10^8$ or $10^9$, as explained below, rather than $P=1699$ that we
achieved for $k=13$, thus requiring roughly $10^5$ to $10^6$ CPU years.

For smaller values of $k$, because of the lower complexity, we were able to
achieve much larger values of $P$, and thus it made sense to set Digits larger.
For example we achieved $P=942,939,827$ for $k=4$ with Digits$=25$, and
$P=1,212,569$ for $k=8$ with Digits$=16$. Our specific choice of Digits and $P$
is listed for each $k$ in the captions of
Tables~\ref{table:coeff_4_5_6}-\ref{table:coeff_13}.

Second, the machine that we used is a
multi-user machine and, even though it has 192 GB RAM, we were running
processes for several $k\leq 13$ simultaneously and also competing for the machine's
resources with other researchers and projects, so we had to temper our use
of the machine's memory and processors.

Note that the memory requirements are substantial, using
\begin{equation}
    O\left( {2k \choose k} k^2 \Digits \right)
\end{equation}
bits for the storage as we loop through $p$ to compute the terms of the terms
in~\eqref{eq:c_r as limit}.

\subsection{Complexity analysis}
\label{sec:complexity}

A rough estimate of the complexity involved in our computations can now be described.
To numerically compute each $c_r$, with $0 \leq r \leq k^2$, we needed to sum
the $2k \choose k$ terms of~\eqref{eq:c_r as limit}. Notice that the dependence on $r$
only appears in the first factor of $H_r$ in~\eqref{eq:H_r}.

The bulk of our computing time was spent in evaluating the $2k \choose k$ Euler products
$A_k(\delta_{\sigma(1)},\ldots,\delta_{\sigma(2k)})$, and we limit ourselves to describing
the complexity of computing the local factors, for all $p\leq P$, in~\eqref{eq:A_k euler product 2},
and the cubic accelerant~\eqref{eq:cubic acc}, to $\approx k^2 \Digits$ decimal places.
The choice of $P$ will be described in the next subsection.

To compute the local factor at $p$ expressed in~\eqref{eq:A_k euler product 2} requires, to begin with,
one exponentiation, namely $p^\delta$, one division to get $p^{-\delta}$, and $O(k)$ multiplications
to compute all relevant powers of $p^{\delta}$ which we then store. In fact, these values can be used,
for given $p$, across all $2k \choose k$ terms of~\eqref{eq:c_r as limit}, and thus forms an
insignificant portion of the overall computation.

We can rewrite the local factor of~\eqref{eq:A_k euler product 2} as
\begin{equation}
    \label{eq:local factor for computing}
    \left(\prod_{i=1}^k \prod_{m=1}^k (1-p^{-1+z_{i+k}-z_m})\right)
    \sum_{j=1}^k
    \left(\prod_{i\neq j} (1-p^{z_{i+k}-z_{j+k}})\right)^{-1}
    \left(\prod_{m=1}^k (1-p^{-1+z_{j+k}-z_m})\right)^{-1}.
\end{equation}
In our application, we set $z_j=\sigma(j) \delta$, with $\sigma \in \Xi$.

The double product outside the sum can be evaluated using our table of values of
powers of $p^\delta$, and $O(k^2)$ multiplications. The number of multiplications
can be reduced to $O(k \log{k})$ if we loop through the indices to count
multiplicities amongst the terms that appear, and use repeated squaring, say,
to carry out the powers of each factor (when counted according to multiplicity).
Recall that there are only $O(k)$ possible values for the factors in the product.

Each term of the above sum over $j$ can be computed using $O(k)$
multiplications and one division. For each $p$, we also pre-computed, stored,
and recycled the possible values, of which there are $O(k)$, of the individual factors that
appear in the sum, as they occur frequently over the $2k \choose k$ terms
in~\eqref{eq:c_r as limit}.
Thus, the sum over $j$ involves, altogether,
$O(k^2)$ multiplications of precomputed values.
The precomputed values were then discarded once we moved onto the next $p$.

Hence, evaluating~\eqref{eq:A_k euler product 2} for all $p\leq P$ can be achieved using
\begin{equation}
    \label{eq:complexity p<=P}
    \ll \frac{P k^2}{\log{P}}
\end{equation}
arithmetic operations (mainly multiplications) on numbers of $O(k^2 \Digits)$ bits.
The factor $P/\log{P}$ is to account, asymptotically, for the number of primes up to $P$.
Using a Fast Fourier Transform method (such as Sch\"onhage-Strassen's algorithm) for
the high precision multiplications, the contribution from
computing the local factors of~\eqref{eq:A_k euler product 2}
for all $p\leq P$ across all $2k \choose k$
terms of~\eqref{eq:c_r as limit} involves, for any $\epsilon>0$,
\begin{eqnarray}
    \label{eq:bit estimate}
    &&\ll {2k \choose k} k^{4+\epsilon} \Digits^{1+\epsilon}\frac{P}{\log{P}} \notag \\
    &&\ll 2^{2k} k^{7/2+\epsilon} \Digits^{1+\epsilon} \frac{P}{\log{P}}
\end{eqnarray}
bit operations, the latter bound following from Stirling's formula applied to
the binomial coefficient.

The evaluation of the rhs of~\eqref{eq:cubic acc} can be carried out using similar methods.
The bulk of the computation is spent in evaluating the product over $p\leq P$
that appears in the denominator. While the product involves $O(k^6)$ factors
(in the product over $i_1,i_2,i_3,j_1,j_2,j_3$) there are, for each $p$, just
$O(k)$ possible values for each factor when $z_j = \sigma(j) \delta$, $1 \leq
\sigma(j) \leq 2k$. Hence these factors appear with high multiplicity.
Furthermore, the multiplicities are independent of $p$, and we can count these
and store the multiplicities for each of the $2k \choose k$ permutations $\sigma$.
The bit complexity in counting these powers is
\begin{eqnarray}
    \ll {2k \choose k} k^{6+\epsilon}.
\end{eqnarray}
Looking ahead to our overall bit complexity in~\eqref{eq:bit estimate b},
we find that counting powers is a small portion of the overall time needed.

Hence, by counting these multiplicities and using repeated squaring, we can
compute the denominator of~\eqref{eq:cubic acc} with $O(k\log{k})$, rather than
$O(k^6)$, high precision multiplications. Thus, the complexity in
computing~\eqref{eq:cubic acc}, for given $P$, is dominated by~\eqref{eq:bit
estimate}.

\subsection{Choice of $P$}
\label{sec:choice P}

We now give an estimate, with an explicit dependence on $k$, in the error in $A_k$ from
truncating its Euler product~\eqref{eq:A_k} at given $P$. The local factor at $p$
of $A_k(\delta_{\sigma(1)},\ldots,\delta_{\sigma(2k)})$ can be expanded in a
series whose terms are of the form $1/p^{r+m\delta}$, where $r\geq 0$ and
$-2rk\delta \leq m \leq 2rk\delta$ (the latter because $\delta_j = j \delta$ with
$1\leq j \leq 2k$).

These arise in~\eqref{eq:uw integral} by matching up an equal number of
$e(\theta)$'s with $e(-\theta)$'s in the integrand of that formula. In the same
way that we worked out the terms up to degree six in~\eqref{eq:uw integral}, we may
consider the terms of degree $2r$ by passing to~\eqref{eq:uw function},
and focusing on the terms that have the same number, $r$, of $u$'s and $w$'s.

We can get an upper bound on the sum total of the terms of degree $2r$ in~\eqref{eq:uw integral}
by examining
the coefficient of $x^{2r}$ in:
\begin{equation}
    \label{eq:uw function c}
    (1+x^2)^{k^2}
    (1-x)^{-2k},
\end{equation}
with $x=p^{-1/2+2k\delta}$. To get this expression,
we have replaced the minus sign in~\eqref{eq:uw function} with a plus sign,
so that both $(1+x^2)^{k^2}$ and $(1-x)^{-2k}$ have Macluarin
series with positive coefficients, and all the $u_j$'s and $w_j$'s with $x$'s.

Now, expanding $(1+x^2)^{k^2}$ and $(1-x)^{-2k}$ in their Maclaurin series,
and multiplying out gives
\begin{equation}
    \label{eq:uw function d}
    (1+x^2)^{k^2}
    (1-x)^{-2k} = \sum_{0}^\infty h_m(k) x^m,
\end{equation}
with
\begin{equation}
    \label{eq:h_m}
    h_m(k) = \sum_{0\leq j \leq m/2} {k^2 \choose j} \frac{2k(2k+1)\ldots(2k+m-2j-1)}{(m-2j)!}.
\end{equation}
Now, because $(2k+l)/l\leq k+1$ if $l\geq 2$, and because $\sum 1/j!$ converges, we have
\begin{equation}
    \label{eq:h_m bound}
    |h_m(k)| \ll  (k+1)^m.
\end{equation}
Our application of this bound is to the case $m=2r$.

In a similar way, we can bound the terms involving $p^{-r}$ on expanding out the products
in~\eqref{eq:prod}, and prove that
\begin{eqnarray}
    \label{eq:asympt local factor}
    &&\prod_p
    \prod_{l,j=1}^k (1-p^{-1-z_l+z_{k+j}})
    \int_0^1 \prod_{j=1}^k
    \left(1-\frac{e(\theta)}{p^{\frac12 +z_j}}\right)^{-1}
    \left(1-\frac{e(-\theta)}{p^{\frac12 -z_{k+j}}}\right)^{-1}\,d\theta \notag \\
    &&= \text{(equation~\eqref{eq:prod})}(1+O(k^8 p^{-4+8k\delta}))
\end{eqnarray}

Hence, the relative error in using~\eqref{eq:cubic acc} to approximate
$A_k(\delta_{\sigma(1)},\ldots,\delta_{\sigma(2k)})$ is
\begin{equation}
    \label{eq:acc error}
    \ll k^8 \sum_{p>P} p^{-4+8k\delta} \ll k^8 \log(P)^{-1} P^{-3+8k\delta}/(3-8k\delta),
\end{equation}
assuming that $3>8k\delta$. Recall that $\delta=10^\Digits$.
Thus, for
\begin{equation}
    P \gg k^{8/3} 10^{\Digits/3},
\end{equation}
and also assuming that
$\Digits$ satisfies, say, $k(\log(k)+\Digits) 10^{-\Digits}< 1/100$
(this easily satisfied condition allows us to ignore the $8k\delta$), we have that
the relative error in $A_k(\delta_{\sigma(1)},\ldots,\delta_{\sigma(2k)})$
from truncating the Euler product at $P$ is $< 10^{-\Digits}$.

Without using~\eqref{eq:cubic acc} to accelerate the convergence of
the Euler product, we would need $P \approx k^4 10^\Digits$ to achieve comparable accuracy,
and using a degree 2 approximation, rather than degree 3, would require
$P \approx k^3 10^{\Digits/2}$. The powers on $k$ here arise as
in~\eqref{eq:asympt local factor}.

Combining~\eqref{eq:bit estimate} with this choice of $P$ gives a bound for
the number of bit operation used
\begin{eqnarray}
    \label{eq:bit estimate b}
    \ll 2^{2k} k^{37/6+\epsilon} \Digits^{\epsilon} 10^{\Digits/3}.
\end{eqnarray}

Note that using a specific desired precision, Digits, does not necessarily
result in that precision being achieved as one also needs to take into account
the implied constant in the $O$ term in~\eqref{eq:c_r as limit b} which depends
both on $k$ and on $r$. Ideally, in place of this $O$-term, one would like to
have uniform bounds with the dependence on $k$ and $r$ explicitly given.
However, this is daunting for several reasons. First, the terms being summed
in~\eqref{eq:c_r as limit b} have poles, with respect to $\delta$, of order
$r$, that cancel, with the most challenging case being when $r=k^2$. To see
one's way through this enormous cancellation involves examining the terms of
the multi-variate series expansions of the summand up to degree $r$. In fact,
our algorithm circumvents this difficulty (experimentally) by using high
precision to account for the high amount of cancellation resulting from the
high order poles that annihilate one another.

One could give extremely crude estimates for this $O$-term by using the
$2k$-fold residue that results from taking, in~\eqref{eq:P_k alpha}, the $r$-th
term on expanding the $\exp$ in its Maclaurin series (with respect to the $x$
variable), and then differentiating, with respect to
the shifts, under the integral sign.

The $2k$-fold residue has the advantage of being analytic when the shifts
$\alpha_j$ are set to 0, i.e. it encodes the cancelling of the poles
in~\eqref{eq:c_r as limit b}. However, this advantage is deceptive, as
the high multiplicity of the factors in the denominator of~\eqref{eq:P_k alpha}
introduces terms of high degree into the analysis. Furthermore, the
$c_r(k)$'s are, for the large part, very small, as can be seen in our tables.
Estimates on the $c_r(k)$'s alone have only been successfully carried out for a
relatively small range of $r$, specifically,  $r< k^\beta$, where
$\beta<1$~\cite{HR}. Presumably the implied constants in the $O$-term are also
comparatively small, and crudely bounding the integrand of the $2k$-fold
residue will thus not produce estimates useful in asserting the numerical
correctness of the coefficients in our tables.


In practice, rather than working with explicit constants in the truncation
bounds, both for the $O$ term in~\eqref{eq:c_r as limit b} and the Euler
product, we experimented by
taking different values of $\delta=10^{-\Digits}$ and $P$, using our
estimates as guides. We inspected the stability of our
numerical values of $c_r(k)$ by comparing those computed for
a given $P$ against those with $P$ replaced by the first prime smaller
than $P/3$, and only outputting the digits that agreed. It seems, from our
tables, that the coefficients $c_r(k)$ with mid-range values of $r$ are more
stable and converge faster with respect to $P$, especially for larger values of $k$.
We did not explore the reason for this, but presumably the lower terms,
beyond those resulting from our cubic accelerant, have comparatively smaller coefficients
for those values of $r$ for which $c_r(k)$ converges faster.

Numerical values of the coefficients for $4 \leq k \leq 13$ thus obtained
are presented in the Tables~\ref{table:coeff_4_5_6}-~\ref{table:coeff_13} below
in scientific `e' notation, for example $1.2e-3 = 1.2\times10^{-3}$.
High precision values of $c_r(k)$ for $k=1,2,3$, can be obtained from~\cite{CFKRS}.

It is also worth mentioning that all the digits of the coefficients computed in
this manner agree, except in a few instances where the last decimal place
differs slightly,
with the results of the first method of~\cite{CFKRS2} (see
Section~\ref{sec:method 1}). That method has the advantage of producing
high precision values of the coefficients, but is limited to relatively
small values of $r$. We reran the program used in~\cite{CFKRS2} for
$k\leq 13$ and $r \leq 7$ and display those values in
Table~\ref{table:method1} for comparison.

\begin{table}
\centerline{
\begin{tabular}{|c|c|c|c|}
\hline
$r$  &$c_r(4)$ & $c_r(5)$ & $c_r(6)$ \cr
\hline
0 & 2.465018391934227354079894e-13 &	1.416001020622731200955e-24 &	5.12947340914919112e-40 \cr
1 & 5.450140573117186559363058e-11 &	7.380412756494451305968e-22 &	5.306732809926444246e-37 \cr
2 & 5.28772963479120311384897e-09 &	1.7797796235196529053094e-19 &	2.6079207711483512396e-34 \cr
3 & 2.96411431799939794596918e-07 &	2.635886609660724758286e-17 &	8.1016132157790177281e-32 \cr
4 & 1.064595006812847051321182e-05 &	2.6840545349997485760134e-15 &	1.7861297380093099773e-29 \cr
5 & 2.5702983342426340235494e-04 &	1.993641309249897180312e-13 &	2.9743167108636063482e-27 \cr
6 & 4.2639216163116947218762e-03 &	1.1184855124933629437778e-11 &	3.8877082911558678876e-25 \cr
7 & 4.89414245142160102712761e-02 &	4.842797553044804165519e-10 &	4.09224261406862935514e-23 \cr
8 & 3.878526654019553499833e-01 &	1.639801308496156099797e-08 &	3.5314663856570325725e-21 \cr
9 & 2.10913382864873355204e+00 &	4.374935105492246330412e-07 &	2.530637690060973478289e-19 \cr
10 & 7.8325356118822623579303e+00 &	9.22633350296530326337e-06 &	1.5198191029685924995e-17 \cr
11 & 1.982806812499890923e+01 &	1.537677778207107946991e-04 &	7.70015137609237458270e-16 \cr
12 & 3.388893203738368856e+01 &	2.01902775807813195907e-03 &	3.3061210414107436046e-14 \cr
13 & 3.82033062189019517e+01 &	2.07727067284846475474e-02 &	1.2064041518984715612e-12 \cr
14 & 2.560441501227035e+01 &	1.6625058643910393652e-01 &	3.7467192541626917996e-11 \cr
15 & 1.06189693794016e+01 &	1.026466777849473756e+00 &	9.9056942856889097902e-10 \cr
16 & 7.089464552244e-01 &	4.848589278343642478e+00 &	2.2273885767179683823e-08 \cr
17 & &	1.73908760901310234e+01 &	4.251372866816786076e-07 \cr
18 & &	4.7040877087561734e+01 &	6.8674335769870947550e-06 \cr
19 & &	9.511661794587886e+01 &	9.351583018775044262e-05 \cr
20 & &	1.41444460064317e+02 &	1.068316421173022528e-03 \cr
21 & &	1.4935694999630e+02 &	1.01807023862361485e-02 \cr
22 & &	1.0588728028422e+02 &	8.04186793058379244e-02 \cr
23 & &	4.41362307288e+01 &	5.2296141941724947e-01 \cr
24 & &	2.010650046e+01 &	2.7802017665195719e+00 \cr
25 & &	-1.2701703e+00 &	1.200111408801811e+01 \cr
26 & & &		4.179670936891264e+01 \cr
27 & & &		1.16723095829484e+02 \cr
28 & & &	2.5939897299715e+02 \cr
29 & & &		4.524908135220e+02 \cr
30 & & &		6.0117334836510e+02 \cr
31 & & &		5.7354384553122e+02 \cr
32 & & &		3.75018676133e+02 \cr
33 & & &		2.46890415605e+02 \cr
34 & & &		2.454954369e+02 \cr
35 & & &		1.603303769e+02 \cr
36 & & &		-3.78219665e+01 \cr
\hline
\end{tabular}
}
\caption{Coefficients $c_r(k)$ for $k = 4,5,6$. For $k=4,5$ we used $\delta=10^{-25}$, and $P=942939827,180343651$
respectively. For $k=6$ we combined two data sets using $\delta=10^{-25}$, $P=25501199$ for
$0 \leq r \leq 29$, and $\delta=10^{-16}$, $P=608121859$ for $30 \leq r \leq 36$.
}
\label{table:coeff_4_5_6}
\end{table}

\begin{table}
\centerline{
\begin{tabular}{|c|c||c|c||c|c|}
\hline
$r$  & $c_r(7)$ & $r$ & $c_r(7)$ & $r$ & $c_r(7)$ \cr
\hline
0 & 6.5822847876005500e-60 & 1 & 1.2041430555451870e-56 & 2 & 1.0621355717492720e-53 \cr
3 & 6.0172653760159300e-51 & 4 & 2.4606287673240130e-48 & 5 & 7.7390121665211530e-46 \cr
6 & 1.9478649494952360e-43 & 7 & 4.0307684926363700e-41 & 8 & 6.9917763337237880e-39 \cr
9 & 1.0314019779812270e-36 & 10 & 1.3082869144993580e-34 & 11 & 1.4392681201435320e-32 \cr
12 & 1.3825312154986080e-30 & 13 & 1.1657759371318020e-28 & 14 & 8.6652476933527220e-27 \cr
15 & 5.6962227424753780e-25 & 16 & 3.3197648540507990e-23 & 17 & 1.7183970393294220e-21 \cr
18 & 7.9096788893235440e-20 & 19 & 3.2396929335740840e-18 & 20 & 1.1809579273268370e-16 \cr
21 & 3.8302270051510130e-15 & 22 & 1.1044706290361260e-13 & 23 & 2.8282258231583490e-12 \cr
24 & 6.4210662257609940e-11 & 25 & 1.2898755567219640e-09 & 26 & 2.2869667400876520e-08 \cr
27 & 3.5683995004969530e-07 & 28 & 4.8834071041615500e-06 & 29 & 5.8391045220798220e-05 \cr
30 & 6.0742037327532430e-04 & 31 & 5.4716438254364890e-03 & 32 & 4.2465903403750590e-02 \cr
33 & 2.824549346789606e-01 & 34 & 1.601333066518585e+00 & 35 & 7.696699614092694e+00 \cr
36 & 3.12035202072577e+01 & 37 & 1.06197143546798e+02 & 38 & 3.019136554174e+02 \cr
39 & 7.117410357280e+02 & 40 & 1.37009445510e+03 & 41 & 2.0827987442e+03 \cr
42 & 2.357363536e+03 & 43 & 1.93463843e+03 & 44 & 1.75714310e+03 \cr
45 & 2.853378e+03 & 46 & 3.100593e+03 & 47 & 3.3940e+02 \cr
48 & -1.20854e+03 & 49 & -5.0194e+02 & & \cr
\hline
\end{tabular}
}
\caption{Coefficients for $k = 7$ truncating at $P = 11015647$, and using $\delta=10^{-16}$.
}
\label{table:coeff_7}
\end{table}

\begin{table}
\centerline{
\begin{tabular}{|c|c||c|c||c|c|}
\hline
$r$ & $c_r(8)$ & $r$ & $c_r(8)$ & $r$ & $c_r(8)$ \cr
\hline
0 & 1.870442160117e-84 & 1 & 5.570219365179e-81 & 2 & 8.0727983790767e-78 \cr
3 & 7.5876025208718e-75 & 4 & 5.2002464291967e-72 & 5 & 2.7705098412043e-69 \cr
6 & 1.1944832708049e-66 & 7 & 4.28399526987474e-64 & 8 & 1.30388783552972e-61 \cr
9 & 3.41901547588078e-59 & 10 & 7.814883564309408e-57 & 11 & 1.5716390984187150e-54 \cr
12 & 2.8019878105779590e-52 & 13 & 4.455999566980015e-50 & 14 & 6.353422665784412e-48 \cr
15 & 8.1564318452061200e-46 & 16 & 9.461598292152723e-44 & 17 & 9.946939863498895e-42 \cr
18 & 9.5005755608004230e-40 & 19 & 8.2610748506972050e-38 & 20 & 6.5505802739681480e-36 \cr
21 & 4.7431624904297110e-34 & 22 & 3.1395134797176380e-32 & 23 & 1.9011136618426250e-30 \cr
24 & 1.0537646192031000e-28 & 25 & 5.3481979466238930e-27 & 26 & 2.4856705239072160e-25 \cr
27 & 1.0578033941072590e-23 & 28 & 4.1205620801239040e-22 & 29 & 1.4685071632966620e-20 \cr
30 & 4.7847077672871950e-19 & 31 & 1.4239651481852300e-17 & 32 & 3.8665596212165030e-16 \cr
33 & 9.566613928613862e-15 & 34 & 2.1534462266724710e-13 & 35 & 4.4024065808181710e-12 \cr
36 & 8.1576484690923080e-11 & 37 & 1.367077736688555e-09 & 38 & 2.0668154599686450e-08 \cr
39 & 2.811291057443101e-07 & 40 & 3.430097716487554e-06 & 41 & 3.741873571285350e-05 \cr
42 & 3.63684791980377e-04 & 43 & 3.137462556407600e-03 & 44 & 2.39287182393225e-02 \cr
45 & 1.60675070999268e-01 & 46 & 9.4588131743684e-01 & 47 & 4.8616355021107e+00 \cr
48 & 2.173103986560e+01 & 49 & 8.417133459453e+01 & 50 & 2.81517268211e+02 \cr
51 & 8.0929177383e+02 & 52 & 1.9821841216e+03 & 53 & 4.05873574e+03 \cr
54 & 6.69566487e+03 & 55 & 8.4203977e+03 & 56 & 8.096360e+03 \cr
57 & 9.4961243e+03 & 58 & 1.99106e+04 & 59 & 3.09087e+04 \cr
60 & 1.3133e+04 & 61 & -2.964e+04 & 62 & -4.0582e+04 \cr
63 & -8.56e+03 & 64 & 4.56e+03 &  & \cr
\hline
\end{tabular}
}
\caption{Coefficients for $k = 8$ truncating at $P = 1212569$, and using $\delta=10^{-16}$.
}
\label{table:coeff_8}
\end{table}

\begin{table}
\centerline{
\begin{tabular}{|c|c||c|c||c|c|}
\hline
$r$ & $c_r(9)$ & $r$ & $c_r(9)$ & $r$ & $c_r(9)$ \cr
\hline
0 & 7.920155238e-114 & 1 & 3.608743873e-110 & 2 & 8.051296272e-107 \cr
3 & 1.172362406e-103 & 4 & 1.2530058769e-100 & 5 & 1.0481427376e-97 \cr
6 & 7.1456310032e-95 & 7 & 4.0821744596e-92 & 8 & 1.9941925290e-89 \cr
9 & 8.4594205088e-87 & 10 & 3.1538144173e-84 & 11 & 1.0433789942e-81 \cr
12 & 3.08731306439e-79 & 13 & 8.22414128961e-77 & 14 & 1.98314875716e-74 \cr
15 & 4.34906593630e-72 & 16 & 8.70843432097e-70 & 17 & 1.597618788502e-67 \cr
18 & 2.693264006511e-65 & 19 & 4.1828312947266e-63 & 20 & 5.9981084135426e-61 \cr
21 & 7.9570365567690e-59 & 22 & 9.7816383607331e-57 & 23 & 1.1159063652824e-54 \cr
24 & 1.1828936996292e-52 & 25 & 1.16636205776320e-50 & 26 & 1.0707494871812e-48 \cr
27 & 9.1588630961799e-47 & 28 & 7.30410594473147e-45 & 29 & 5.43352058419794e-43 \cr
30 & 3.77180308962311e-41 & 31 & 2.44391067332666e-39 & 32 & 1.47828492081742e-37 \cr
33 & 8.34814115852835e-36 & 34 & 4.401070135763961e-34 & 35 & 2.16570032669324e-32 \cr
36 & 9.94493396339446e-31 & 37 & 4.26008480938235e-29 & 38 & 1.7015865979528490e-27 \cr
39 & 6.33392559975936e-26 & 40 & 2.19581041600330e-24 & 41 & 7.08426675189073e-23 \cr
42 & 2.1252127802020740e-21 & 43 & 5.92241988702831e-20 & 44 & 1.53150283306451e-18 \cr
45 & 3.670620764266784e-17 & 46 & 8.14314143690125e-16 & 47 & 1.66973151449710e-14 \cr
48 & 3.15948278435009e-13 & 49 & 5.507449862131347e-12 & 50 & 8.82746086945955e-11 \cr
51 & 1.29834184462004e-09 & 52 & 1.74846308044295e-08 & 53 & 2.1508513613974e-07 \cr
54 & 2.4107241046116e-06 & 55 & 2.4552071751795e-05 & 56 & 2.265572158192e-04 \cr
57 & 1.888394621727e-03 & 58 & 1.4172609520338e-02 & 59 & 9.546112504952e-02 \cr
60 & 5.751640191248e-01 & 61 & 3.08994036050e+00 & 62 & 1.47569216973e+01 \cr
63 & 6.2480977623e+01 & 64 & 2.3393469844e+02 & 65 & 7.721276404e+02 \cr
66 & 2.23430645e+03 & 67 & 5.60365351e+03 & 68 & 1.1907448e+04 \cr
69 & 2.062257e+04 & 70 & 2.776900e+04 & 71 & 3.06185e+04 \cr
72 & 4.7187e+04 & 73 & 1.2202e+05 & 74 & 2.314e+05 \cr
75 & 1.255e+05 & 76 & -4.658e+05 & 77 & -1.07e+06 \cr
78 & -5.794e+05 & 79 & 6.75e+05 & 80 & 8.27e+05 \cr
81 & 1.3e+05 & & & & \cr
\hline
\end{tabular}
}
\caption{Coefficients for $k = 9$ truncating at $P = 170741$, and using $\delta=10^{-16}$.
}
\label{table:coeff_9}
\end{table}

\begin{table}
\centerline{
\begin{tabular}{|c|c||c|c||c|c|}
\hline
$r$ & $c_r(10)$ & $r$ & $c_r(10)$ & $r$ & $c_r(10)$ \cr
\hline
0 & 3.54888492477e-148 & 1 & 2.35769133101e-144 & 2 & 7.70233663026e-141 \cr
3 & 1.64948634407e-137 & 4 & 2.60451944693e-134 & 5 & 3.23366677841e-131 \cr
6 & 3.28765141574e-128 & 7 & 2.81472946999e-125 & 8 & 2.07111222708e-122 \cr
9 & 1.330223430450e-119 & 10 & 7.54902089850e-117 & 11 & 3.822610704580e-114 \cr
12 & 1.741117024160e-111 & 13 & 7.181397221310e-109 & 14 & 2.697527199940e-106 \cr
15 & 9.272621188300e-104 & 16 & 2.929091296840e-101 & 17 & 8.533614631040e-99 \cr
18 & 2.300282495690e-96 & 19 & 5.752933603890e-94 & 20 & 1.338222245340e-91 \cr
21 & 2.901664262540e-89 & 22 & 5.876116434230e-87 & 23 & 1.113295629950e-84 \cr
24 & 1.976420578280e-82 & 25 & 3.292284406610e-80 & 26 & 5.152298303310e-78 \cr
27 & 7.583498206660e-76 & 28 & 1.050825438590e-73 & 29 & 1.372032297570e-71 \cr
30 & 1.689305278750e-69 & 31 & 1.962725749490e-67 & 32 & 2.153176878340e-65 \cr
33 & 2.231491144780e-63 & 34 & 2.185751487190e-61 & 35 & 2.024236345360e-59 \cr
36 & 1.773015880790e-57 & 37 & 1.469142759970e-55 & 38 & 1.151856372240e-53 \cr
39 & 8.546210010780e-52 & 40 & 6.000996648940e-50 & 41 & 3.988024973250e-48 \cr
42 & 2.508210346100e-46 & 43 & 1.492814370840e-44 & 44 & 8.406690331770e-43 \cr
45 & 4.478585581910e-41 & 46 & 2.256576452310e-39 & 47 & 1.075045360230e-37 \cr
48 & 4.840860160570e-36 & 49 & 2.059518250610e-34 & 50 & 8.274865115540e-33 \cr
51 & 3.138261284280e-31 & 52 & 1.122809702480e-29 & 53 & 3.787412159890e-28 \cr
54 & 1.203656074660e-26 & 55 & 3.601325369860e-25 & 56 & 1.013609245800e-23 \cr
57 & 2.681299896780e-22 & 58 & 6.659997348320e-21 & 59 & 1.551711921570e-19 \cr
60 & 3.387479938260e-18 & 61 & 6.920780548090e-17 & 62 & 1.321576287680e-15 \cr
63 & 2.355560477620e-14 & 64 & 3.913136454580e-13 & 65 & 6.049282423090e-12 \cr
66 & 8.687675538270e-11 & 67 & 1.157038978680e-09 & 68 & 1.426296610210e-08 \cr
69 & 1.624086810140e-07 & 70 & 1.704560517010e-06 & 71 & 1.645250259970e-05 \cr
72 & 1.456896715700e-04 & 73 & 1.180639233900e-03 & 74 & 8.733236736800e-03 \cr
75 & 5.881138307580e-02 & 76 & 3.596179964360e-01 & 77 & 1.991704581680e+00 \cr
78 & 9.96798553888e+00 & 79 & 4.49881033284e+01 & 80 & 1.8275974779e+02 \cr
81 & 6.6682890321e+02 & 82 & 2.177262632e+03 & 83 & 6.31448075e+03 \cr
84 & 1.60290666e+04 & 85 & 3.47003569e+04 & 86 & 6.163495e+04 \cr
87 & 8.726963e+04 & 88 & 1.147239e+05 & 89 & 2.48873e+05 \cr
90 & 7.40512e+05 & 91 & 1.4296e+06 & 92 & 2.56e+05 \cr
93 & -6.2748e+06 & 94 & -1.489e+07 & 95 & -7.97e+06 \cr
96 & 2.260e+07 & 97 & 4.02e+07 & 98 & 1.43e+07 \cr
99 & -1.08e+07 & 100 & -5.22e+06 &  & \cr
\hline
\end{tabular}
}
\caption{Coefficients for $k = 10$ truncating at $P = 675929$, and using $\delta=10^{-16}$.
}
\label{table:coeff_10}
\end{table}

\begin{table}
\centerline{
\begin{tabular}{|c|c||c|c||c|c||c|c|}
\hline
$r$ & $c_r(11)$ & $r$ & $c_r(11)$ & $r$ & $c_r(11)$ & $r$ & $c_r(11)$ \cr
\hline
0 & 1.2451314e-187 & 1 & 1.16057289e-183 & 2 & 5.33593693e-180 & 3 & 1.61328064e-176 \cr
4 & 3.60796892e-173 & 5 & 6.36556626e-170 & 6 & 9.22782529e-167 & 7 & 1.13037084e-163 \cr
8 & 1.19424156e-160 & 9 & 1.10531689e-157 & 10 & 9.07262368e-155 & 11 & 6.67001789e-152 \cr
12 & 4.42795938e-149 & 13 & 2.67250329e-146 & 14 & 1.47494511e-143 & 15 & 7.48038633e-141 \cr
16 & 3.50123968e-138 & 17 & 1.51807108e-135 & 18 & 6.117325360e-133 & 19 & 2.297695573e-130 \cr
20 & 8.065058896e-128 & 21 & 2.651612616e-125 & 22 & 8.182770981e-123 & 23 & 2.374593250e-120 \cr
24 & 6.490958842e-118 & 25 & 1.673872653e-115 & 26 & 4.077854544e-113 & 27 & 9.396885825e-111 \cr
28 & 2.050582239e-108 & 29 & 4.2419484994e-106 & 30 & 8.3264998645e-104 & 31 & 1.5521939465e-101 \cr
32 & 2.7501740419e-99 & 33 & 4.6346698213e-97 & 34 & 7.4337414866e-95 & 35 & 1.13549625358e-92 \cr
36 & 1.65268057424e-90 & 37 & 2.29313645275e-88 & 38 & 3.03459182593e-86 & 39 & 3.831526640190e-84 \cr
40 & 4.617401084980e-82 & 41 & 5.312669020390e-80 & 42 & 5.837617513170e-78 & 43 & 6.127289274640e-76 \cr
44 & 6.144657440350e-74 & 45 & 5.888372616130e-72 & 46 & 5.392853053390e-70 & 47 & 4.720750019910e-68 \cr
48 & 3.950047958650e-66 & 49 & 3.159433634290e-64 & 50 & 2.415656293490e-62 & 51 & 1.765512476040e-60 \cr
52 & 1.233372385020e-58 & 53 & 8.235138259720e-57 & 54 & 5.254770626530e-55 & 55 & 3.203931860340e-53 \cr
56 & 1.866325220220e-51 & 57 & 1.038439545240e-49 & 58 & 5.517819788010e-48 & 59 & 2.799212525590e-46 \cr
60 & 1.355385595070e-44 & 61 & 6.261980454760e-43 & 62 & 2.759510005390e-41 & 63 & 1.159469298800e-39 \cr
64 & 4.643172007710e-38 & 65 & 1.771352314140e-36 & 66 & 6.434590479390e-35 & 67 & 2.224528434850e-33 \cr
68 & 7.315010837230e-32 & 69 & 2.286618699140e-30 & 70 & 6.790446084110e-29 & 71 & 1.914406176510e-27 \cr
72 & 5.120202491620e-26 & 73 & 1.298144910520e-24 & 74 & 3.117357525500e-23 & 75 & 7.084341566220e-22 \cr
76 & 1.522160095850e-20 & 77 & 3.089177852240e-19 & 78 & 5.915547597220e-18 & 79 & 1.067669213960e-16 \cr
80 & 1.814092977890e-15 & 81 & 2.898172629520e-14 & 82 & 4.347699833250e-13 & 83 & 6.115899813510e-12 \cr
84 & 8.055381828660e-11 & 85 & 9.918802503370e-10 & 86 & 1.139891722860e-08 & 87 & 1.220516286130e-07 \cr
88 & 1.215356440380e-06 & 89 & 1.123334898600e-05 & 90 & 9.618099818200e-05 & 91 & 7.612754647750e-04 \cr
92 & 5.558326140280e-03 & 93 & 3.73566507420e-02 & 94 & 2.3062374960e-01 & 95 & 1.3052451738e+00 \cr
96 & 6.7601261642e+00 & 97 & 3.199030312e+01 & 98 & 1.381287107e+02 & 99 & 5.43353076e+02 \cr
100 & 1.94226636e+03 & 101 & 6.2767681e+03 & 102 & 1.8153278e+04 & 103 & 4.6140024e+04 \cr
104 & 1.001842e+05 & 105 & 1.79901e+05 & 106 & 2.732785e+05 & 107 & 4.8073e+05 \cr
108 & 1.4164e+06 & 109 & 4.18178e+06 & 110 & 6.523e+06 & 111 & -7.31e+06 \cr
112 & -6.295e+07 & 113 & -1.26e+08 & 114 & -1.2e+07 & 115 & 4.21e+08 \cr
116 & 7.32e+08 & 117 & 1.7e+08 & 118 & -8.1e+08 & 119 & -8.2e+08 \cr
120 & -1.1e+08 & 121 & 9.9e+07 & & & &\cr
\hline
\end{tabular}
}
\caption{Coefficients for $k = 11$ truncating at $P = 85889$, and using $\delta=10^{-12}$.
}
\label{table:coeff_11}
\end{table}

\begin{table}
\centerline{
\begin{tabular}{|c|c||c|c||c|c||c|c|}
\hline
$r$ & $c_r(12)$ & $r$ & $c_r(12)$ & $r$ & $c_r(12)$ & $r$ & $c_r(12)$\cr
\hline
0 & 2.61438e-232 & 1 & 3.31314e-228 & 2 & 2.07583e-224 & 3 & 8.57284e-221 \cr
4 & 2.62512e-217 & 5 & 6.35698e-214 & 6 & 1.26799e-210 & 7 & 2.14259e-207 \cr
8 & 3.13061e-204 & 9 & 4.01773e-201 & 10 & 4.58505e-198 & 11 & 4.69936e-195 \cr
12 & 4.36135e-192 & 13 & 3.69038e-189 & 14 & 2.86364e-186 & 15 & 2.04802e-183 \cr
16 & 1.35583e-180 & 17 & 8.34021e-178 & 18 & 4.78307e-175 & 19 & 2.56498e-172 \cr
20 & 1.28962e-169 & 21 & 6.09352e-167 & 22 & 2.71168e-164 & 23 & 1.13871e-161 \cr
24 & 4.52018e-159 & 25 & 1.698887e-156 & 26 & 6.054452e-154 & 27 & 2.048658e-151 \cr
28 & 6.589967e-149 & 29 & 2.017472e-146 & 30 & 5.884306e-144 & 31 & 1.636682e-141 \cr
32 & 4.345107e-139 & 33 & 1.101947e-136 & 34 & 2.671618e-134 & 35 & 6.196518e-132 \cr
36 & 1.375825e-129 & 37 & 2.926060e-127 & 38 & 5.964172e-125 & 39 & 1.165705e-122 \cr
40 & 2.185779e-120 & 41 & 3.9336446e-118 & 42 & 6.7972368e-116 & 43 & 1.1281875e-113 \cr
44 & 1.7992516e-111 & 45 & 2.7580481e-109 & 46 & 4.0647787e-107 & 47 & 5.7611834e-105 \cr
48 & 7.8547080e-103 & 49 & 1.0303518e-100 & 50 & 1.3006541e-98 & 51 & 1.58027693e-96 \cr
52 & 1.84826712e-94 & 53 & 2.08119455e-92 & 54 & 2.25644465e-90 & 55 & 2.35580895e-88 \cr
56 & 2.36859061e-86 & 57 & 2.293494547e-84 & 58 & 2.138841512e-82 & 59 & 1.921057763e-80 \cr
60 & 1.6618175746e-78 & 61 & 1.3845226981e-76 & 62 & 1.110902429e-74 & 63 & 8.583981848e-73 \cr
64 & 6.3871821912e-71 & 65 & 4.5761550941e-69 & 66 & 3.1565961813e-67 & 67 & 2.09610369144e-65 \cr
68 & 1.33974800679e-63 & 69 & 8.2410737927e-62 & 70 & 4.8777660410e-60 & 71 & 2.7774938619e-58 \cr
72 & 1.5212162503e-56 & 73 & 8.01193333444e-55 & 74 & 4.05683148804e-53 & 75 & 1.97436017587e-51 \cr
76 & 9.23280422076e-50 & 77 & 4.14742073168e-48 & 78 & 1.78904091502e-46 & 79 & 7.40817919355e-45 \cr
80 & 2.943700133510e-43 & 81 & 1.122014888380e-41 & 82 & 4.10061125513e-40 & 83 & 1.43633249117e-38 \cr
84 & 4.81968183211e-37 & 85 & 1.548559986510e-35 & 86 & 4.76169394977e-34 & 87 & 1.400502643990e-32 \cr
88 & 3.93774707023e-31 & 89 & 1.05777224840e-29 & 90 & 2.71295150626e-28 & 91 & 6.63907833056e-27 \cr
92 & 1.549117947290e-25 & 93 & 3.44390394051e-24 & 94 & 7.28902479104e-23 & 95 & 1.467519582110e-21 \cr
96 & 2.808147993620e-20 & 97 & 5.10250621963e-19 & 98 & 8.79548303360e-18 & 99 & 1.43685110544e-16 \cr
100 & 2.22218436508e-15 & 101 & 3.24998071837e-14 & 102 & 4.489567612890e-13 & 103 & 5.85079932579e-12 \cr
104 & 7.18370974878e-11 & 105 & 8.29873790933e-10 & 106 & 9.0070339192e-09 & 107 & 9.1707380762e-08 \cr
108 & 8.7457189835e-07 & 109 & 7.7990666196e-06 & 110 & 6.4924284455e-05 & 111 & 5.036525894e-04 \cr
112 & 3.6345151422e-03 & 113 & 2.4355146991e-02 & 114 & 1.5129598801e-01 & 115 & 8.699064420e-01 \cr
116 & 4.62298136e+00 & 117 & 2.268127865e+01 & 118 & 1.02629299e+02 & 119 & 4.2781349e+02 \cr
120 & 1.6399167e+03 & 121 & 5.7590488e+03 & 122 & 1.839118e+04 & 123 & 5.26980e+04 \cr
124 & 1.3266774e+05 & 125 & 2.8583e+05 & 126 & 5.2381e+05 & 127 & 9.3264e+05 \cr
128 & 2.342e+06 & 129 & 7.741e+06 & 130 & 1.944e+07 & 131 & 1.41e+07 \cr
132 & -1.08e+08 & 133 & -4.566e+08 & 134 & -6.14e+08 & 135 & 1.1e+09 \cr
136 & 5.58e+09 & 137 & 7.3e+09 & 138 & -6.6e+09 & 139 & -3.37e+10 \cr
140 & -3.7e+10 & 141 & 5.5e+09 & 142 & 4.3e+10 & 143 & 2.7e+10 \cr
144 & 1.9e+09 & & & & & & \cr
\hline
\end{tabular}
}
\caption{Coefficients for $k = 12$ truncating at $P = 12979$, and using $\delta=10^{-12}$.
}
\label{table:coeff_12}
\end{table}

\begin{table}
\centerline{
\begin{tabular}{|c|c||c|c||c|c||c|c||c|c|}
\hline
$r$ & $c_r(13)$ & $r$ & $c_r(13)$ & $r$ & $c_r(13)$ & $r$ & $c_r(13)$ \cr
\hline
0 & 2.58e-282 & 1 & 4.33e-278 & 2 & 3.60e-274 & 3 & 1.97e-270 \cr
4 & 8.05e-267 & 5 & 2.60e-263 & 6 & 6.94e-260 & 7 & 1.57e-256 \cr
8 & 3.08e-253 & 9 & 5.31e-250 & 10 & 8.16e-247 & 11 & 1.13e-243 \cr
12 & 1.42e-240 & 13 & 1.63e-237 & 14 & 1.71e-234 & 15 & 1.67e-231 \cr
16 & 1.51e-228 & 17 & 1.27e-225 & 18 & 9.96e-223 & 19 & 7.34e-220 \cr
20 & 5.08e-217 & 21 & 3.31e-214 & 22 & 2.04e-211 & 23 & 1.19e-208 \cr
24 & 6.56e-206 & 25 & 3.43e-203 & 26 & 1.71e-200 & 27 & 8.11e-198 \cr
28 & 3.67e-195 & 29 & 1.58e-192 & 30 & 6.52e-190 & 31 & 2.57e-187 \cr
32 & 9.69e-185 & 33 & 3.50e-182 & 34 & 1.21e-179 & 35 & 4.04e-177 \cr
36 & 1.29e-174 & 37 & 3.96e-172 & 38 & 1.17e-169 & 39 & 3.32e-167 \cr
40 & 9.068e-165 & 41 & 2.387e-162 & 42 & 6.054e-160 & 43 & 1.480e-157 \cr
44 & 3.490e-155 & 45 & 7.939e-153 & 46 & 1.743e-150 & 47 & 3.695e-148 \cr
48 & 7.563e-146 & 49 & 1.496e-143 & 50 & 2.858e-141 & 51 & 5.280e-139 \cr
52 & 9.429e-137 & 53 & 1.628e-134 & 54 & 2.720e-132 & 55 & 4.396e-130 \cr
56 & 6.874e-128 & 57 & 1.040e-125 & 58 & 1.524e-123 & 59 & 2.160e-121 \cr
60 & 2.966e-119 & 61 & 3.942e-117 & 62 & 5.074e-115 & 63 & 6.3253e-113 \cr
64 & 7.6365e-111 & 65 & 8.9299e-109 & 66 & 1.0115e-106 & 67 & 1.1098e-104 \cr
68 & 1.1795e-102 & 69 & 1.2144e-100 & 70 & 1.2112e-98 & 71 & 1.1702e-96 \cr
72 & 1.0952e-94 & 73 & 9.9295e-93 & 74 & 8.7200e-91 & 75 & 7.4173e-89 \cr
76 & 6.11074e-87 & 77 & 4.87569e-85 & 78 & 3.76736e-83 & 79 & 2.81878e-81 \cr
80 & 2.04204e-79 & 81 & 1.43219e-77 & 82 & 9.72340e-76 & 83 & 6.38944e-74 \cr
84 & 4.063248e-72 & 85 & 2.500254e-70 & 86 & 1.488419e-68 & 87 & 8.5708030e-67 \cr
88 & 4.7730131e-65 & 89 & 2.5701193e-63 & 90 & 1.3378668e-61 & 91 & 6.7309343e-60 \cr
92 & 3.2721950e-58 & 93 & 1.5367284e-56 & 94 & 6.9700536e-55 & 95 & 3.0523694e-53 \cr
96 & 1.29025276e-51 & 97 & 5.26280865e-50 & 98 & 2.0707478e-48 & 99 & 7.8570340e-47 \cr
100 & 2.8738212e-45 & 101 & 1.01291060e-43 & 102 & 3.43894917e-42 & 103 & 1.124213471e-40 \cr
104 & 3.53717984e-39 & 105 & 1.07067892e-37 & 106 & 3.11641068e-36 & 107 & 8.71833921e-35 \cr
108 & 2.34302639e-33 & 109 & 6.04582505e-32 & 110 & 1.49702508e-30 & 111 & 3.55505913e-29 \cr
112 & 8.09182488e-28 & 113 & 1.76422359e-26 & 114 & 3.68197859e-25 & 115 & 7.35070729e-24 \cr
116 & 1.40275907e-22 & 117 & 2.55690377e-21 & 118 & 4.44814493e-20 & 119 & 7.37932244e-19 \cr
120 & 1.16640624e-17 & 121 & 1.755027226e-16 & 122 & 2.51133489e-15 & 123 & 3.41411712e-14 \cr
124 & 4.40504845e-13 & 125 & 5.38821792e-12 & 126 & 6.24112402e-11 & 127 & 6.83724634e-10 \cr
128 & 7.07544905e-09 & 129 & 6.90733979e-08 & 130 & 6.352697060e-07 & 131 & 5.49643677e-06 \cr
132 & 4.4673281e-05 & 133 & 3.40576382e-04 & 134 & 2.4318473e-03 & 135 & 1.6239805e-02 \cr
136 & 1.0128549e-01 & 137 & 5.8923590e-01 & 138 & 3.194017e+00 & 139 & 1.611797e+01 \cr
140 & 7.566471e+01 & 141 & 3.301742e+02 & 142 & 1.337419e+03 & 143 & 5.014434e+03 \cr
144 & 1.73025e+04 & 145 & 5.43813e+04 & 146 & 1.53198e+05 & 147 & 3.7902e+05 \cr
148 & 8.147e+05 & 149 & 1.6140e+06 & 150 & 3.743e+06 & 151 & 1.18e+07 \cr
152 & 3.577e+07 & 153 & 6.04e+07 & 154 & -8.98e+07 & 155 & -8.78e+08 \cr
156 & -2.1e+09 & 157 & 1.30e+09 & 158 & 2.1e+10 & 159 & 5.2e+10 \cr
160 & -9e+09 & 161 & -3.6e+11 & 162 & -8.7e+11 & 163 & -5.2e+11 \cr
164 & 1.6e+12 & 165 & 3.83e+12 & 166 & 2.6e+12 & 167 & -9e+11 \cr
168 & -1.65e+12 & 169 & -3.7e+11 & & & & \cr
\hline
\end{tabular}
}
\caption{Coefficients for $k = 13$ truncating at $P = 1699$, and using $\delta=10^{-12}$.
}
\label{table:coeff_13}
\end{table}

\begin{table}
\centerline{
\begin{tabular}{|c|c|c||c|c|c|}
\hline
$k$ & $r$ & $c_r(k)$ & $k$ & $r$ & $c_r(k)$ \cr
\hline
4 & 0 & .24650183919342273540799e-12 & 9 & 0 & .79201552383685290316e-113 \cr
4 & 1 & .545014057311718655936e-10 & 9 & 1 & .36087438729455558616e-109 \cr
4 & 2 & .5287729634791203113849e-8 & 9 & 2 & .80512962716760934894e-106 \cr
4 & 3 & .29641143179993979459691e-6 & 9 & 3 & .11723624058166636900e-102 \cr
4 & 4 & .10645950068128470513211e-4 & 9 & 4 & .12530058768923713471e-99 \cr
4 & 5 & .2570298334242634023549e-3 & 9 & 5 & .10481427375523351016e-96 \cr
4 & 6 & .426392161631169472187e-2 & 9 & 6 & .71456310032205157639e-94 \cr
4 & 7 & .4894142451421601027126e-1 & 9 & 7 & .40821744596370636463e-91 \cr
\hline
5 & 0 & .14160010206227312010e-23 & 10 & 0 & .35488849247730348098e-147 \cr
5 & 1 & .73804127564944513060e-21 & 10 & 1 & .23576913310137009644e-143 \cr
5 & 2 & .17797796235196529053e-18 & 10 & 2 & .77023366302575780180e-140 \cr
5 & 3 & .26358866096607247583e-16 & 10 & 3 & .16494863440733411303e-136 \cr
5 & 4 & .26840545349997485760e-14 & 10 & 4 & .26045194469316625626e-133 \cr
5 & 5 & .19936413092498971803e-12 & 10 & 5 & .32336667784065596864e-130 \cr
5 & 6 & .11184855124933629438e-10 & 10 & 6 & .32876514157441589044e-127 \cr
5 & 7 & .48427975530448041655e-9 & 10 & 7 & .28147294699934449064e-124 \cr
\hline
6 & 0 & .51294734091491911243e-39 & 11 & 0 & .124513138816594309e-186 \cr
6 & 1 & .53067328099264442456e-36 & 11 & 1 & .116057289076806867e-182 \cr
6 & 2 & .26079207711483512396e-33 & 11 & 2 & .533593692953085661e-179 \cr
6 & 3 & .81016132157790177281e-31 & 11 & 3 & .161328064239033845e-175 \cr
6 & 4 & .17861297380093099773e-28 & 11 & 4 & .360796891855797563e-172 \cr
6 & 5 & .29743167108636063482e-26 & 11 & 5 & .636556626245602757e-169 \cr
6 & 6 & .38877082911558678876e-24 & 11 & 6 & .922782528634884471e-166 \cr
6 & 7 & .40922426140686293551e-22 & 11 & 7 & .113037084302453487e-162 \cr
\hline
7 & 0 & .65822847876005499378e-59 & 12 & 0 & .26143756530064042e-231 \cr
7 & 1 & .12041430555451865785e-55 & 12 & 1 & .33131354510381580e-227 \cr
7 & 2 & .10621355717492716606e-52 & 12 & 2 & .20758311485643071e-223 \cr
7 & 3 & .60172653760159300486e-50 & 12 & 3 & .85728395653185504e-220 \cr
7 & 4 & .24606287673240130820e-47 & 12 & 4 & .26251165314413802e-216 \cr
7 & 5 & .77390121665211526042e-45 & 12 & 5 & .63569771378458374e-213 \cr
7 & 6 & .19478649494952357456e-42 & 12 & 6 & .12679941162107841e-209 \cr
7 & 7 & .40307684926363697065e-40 & 12 & 7 & .21425932836667245e-206 \cr
\hline
8 & 0 & .18704421601168844202e-83 & 13 & 0 & .2577425553942569e-281 \cr
8 & 1 & .55702193651787573285e-80 & 13 & 1 & .4326313738224894e-277 \cr
8 & 2 & .80727983790767114280e-77 & 13 & 2 & .3596648214485737e-273 \cr
8 & 3 & .75876025208717340817e-74 & 13 & 3 & .1974403388369282e-269 \cr
8 & 4 & .52002464291967325362e-71 & 13 & 4 & .8051097934153342e-266 \cr
8 & 5 & .27705098412043360903e-68 & 13 & 5 & .2601086374395006e-262 \cr
8 & 6 & .11944832708048773794e-65 & 13 & 6 & .6934795975034907e-259 \cr
8 & 7 & .42839952698747409380e-63 & 13 & 7 & .1569263805950391e-255 \cr
\hline
\end{tabular}
}
\caption{For comparison, high precision values of $c_r(k)$, $4 \leq k \leq 13$, $r\leq 7$,
computed using the program for method 1 of~\cite{CFKRS}. The values here agree with those
in Tables~\ref{table:coeff_4_5_6}-~\ref{table:coeff_13} computed using our cubic accelerant
method with experimentally determined remainders,
except they are occasionally slightly off in the last decimal place.
}
\label{table:method1}
\end{table}

\section{\textbf{Checking the Moment Polynomial Conjectures}}
\label{sec:checking conj}

Let
\begin{equation}
    \label{eq:actual moments}
    \text{Data$_k(T)$}= \int_{T_0}^\text{T} |\zeta(1/2+it)|^{2k} dt,
\end{equation}
with $T_0=14.134725\ldots$ being the imaginary part of the first non-trivial zero
of zeta, and let
\begin{equation}
    \label{eq:conjectured moments}
    \text{Conjecture$_k(T)$}= \int_{T_0}^\text{T} P_k(\log \frac{t}{2 \pi}) dt.
\end{equation}
We used our tables of $c_r(k)$ and a simple integration by parts
to compute the latter for given $T$. Note that because the leading coefficients of
$P_k(x)$ are very small, in the range of $T$ considered, the $c_r(k)$'s with midrange and
higher values of $r$ contribute the dominant amount to the integral. To accurately compute the
prediction of~\cite{CFKRS} one does need the lower terms of the polynomial $P_k(x)$.
Table~\ref{table:max} displays, for $1 \leq k \leq 13$ and $T=10^8$,
the $r$ for which the corresponding term in~\eqref{eq:conjectured moments}
contributes the dominant amount to integral.

\begin{figure}[h]
    \centerline{
        \includegraphics[width=.58\textwidth,height=2.7in]{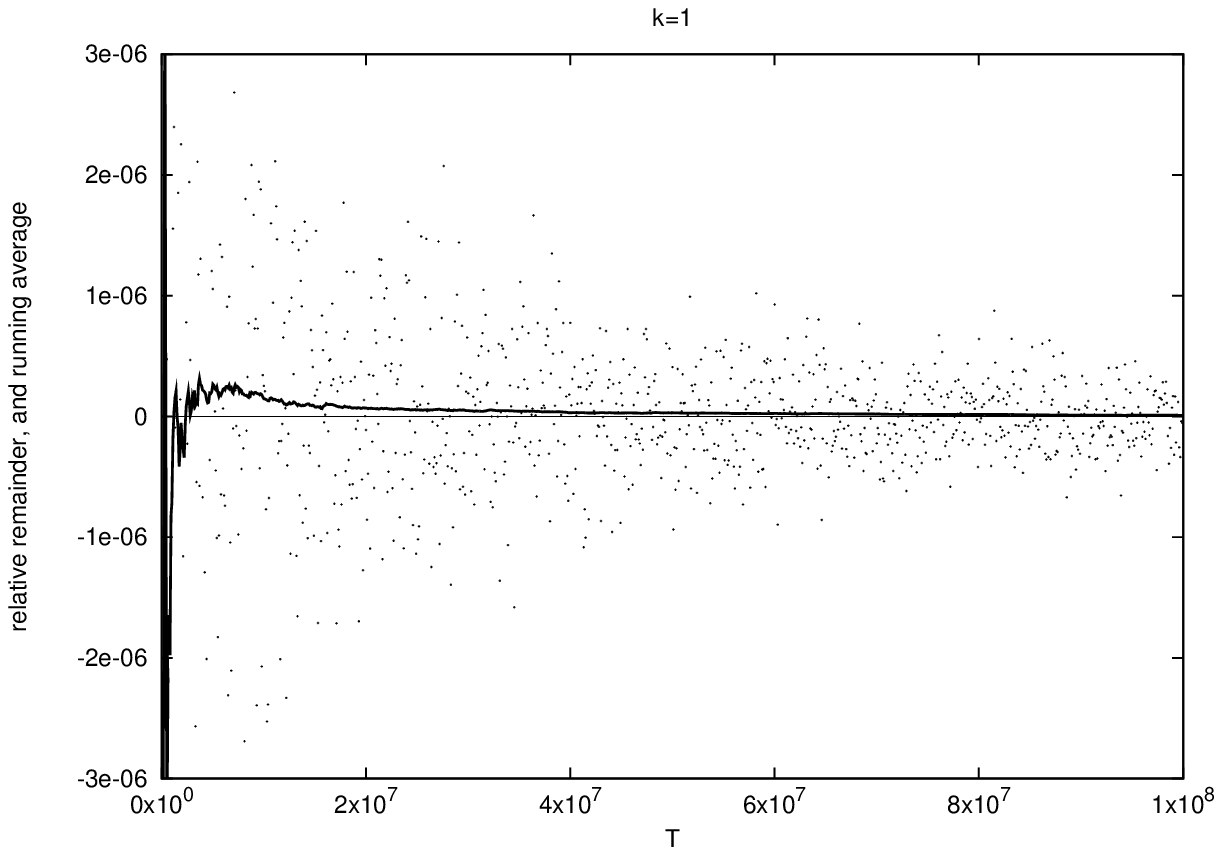}
        \includegraphics[width=.58\textwidth,height=2.7in]{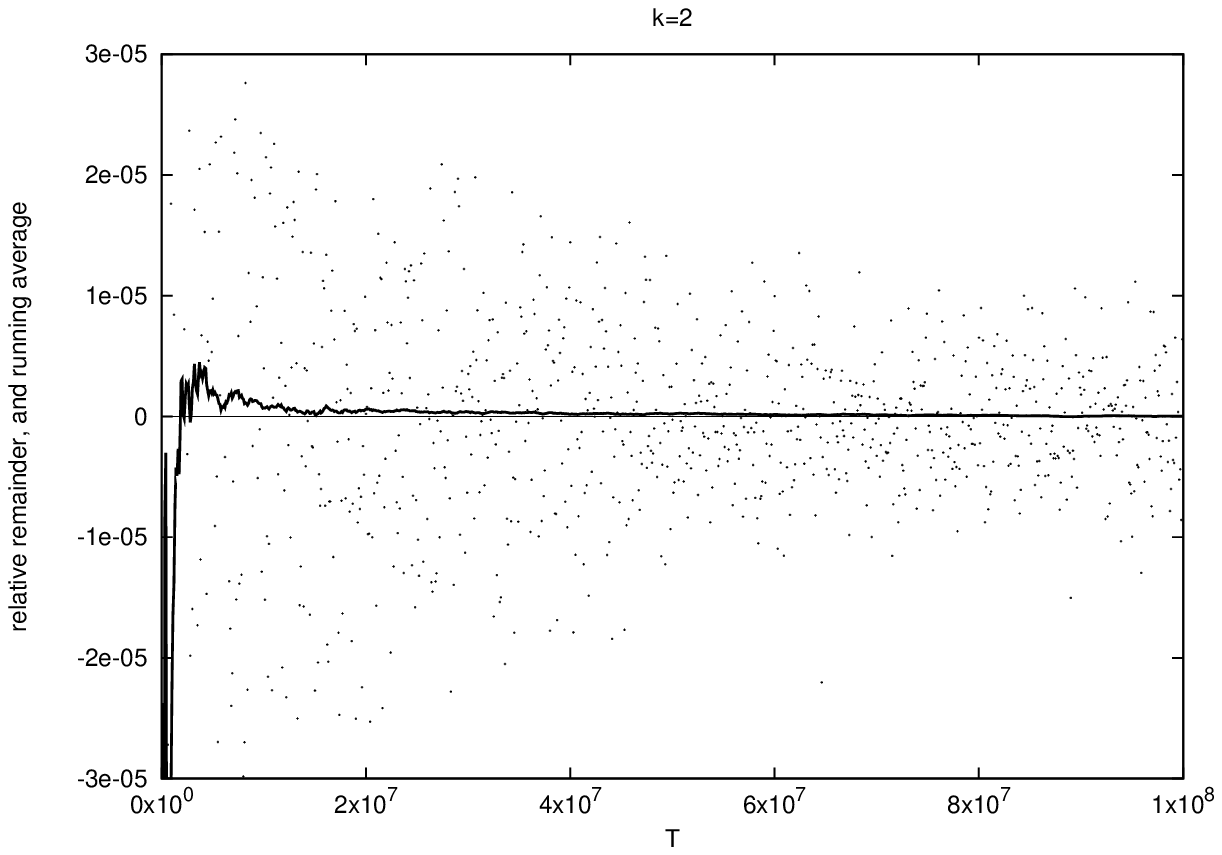}
    }
    \centerline{
       \includegraphics[width=.58\textwidth,height=2.7in]{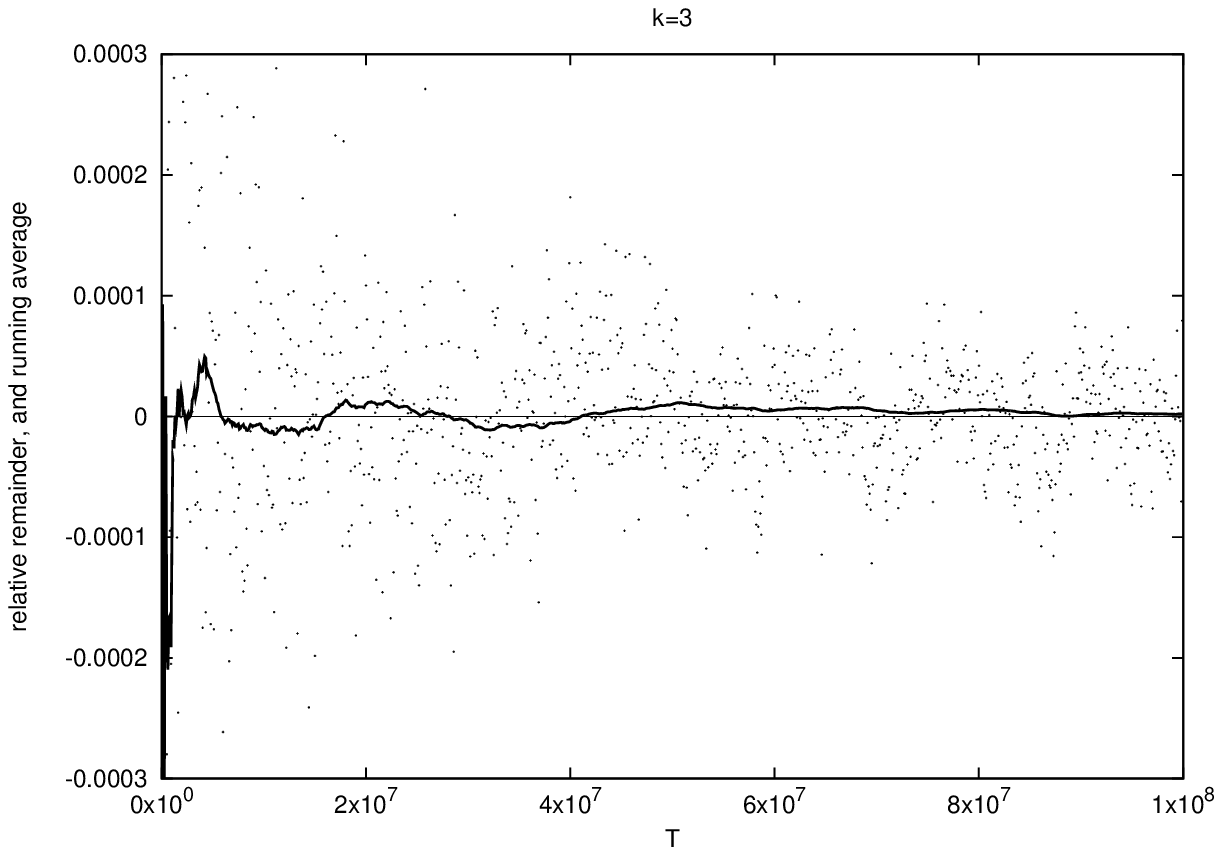}
       \includegraphics[width=.58\textwidth,height=2.7in]{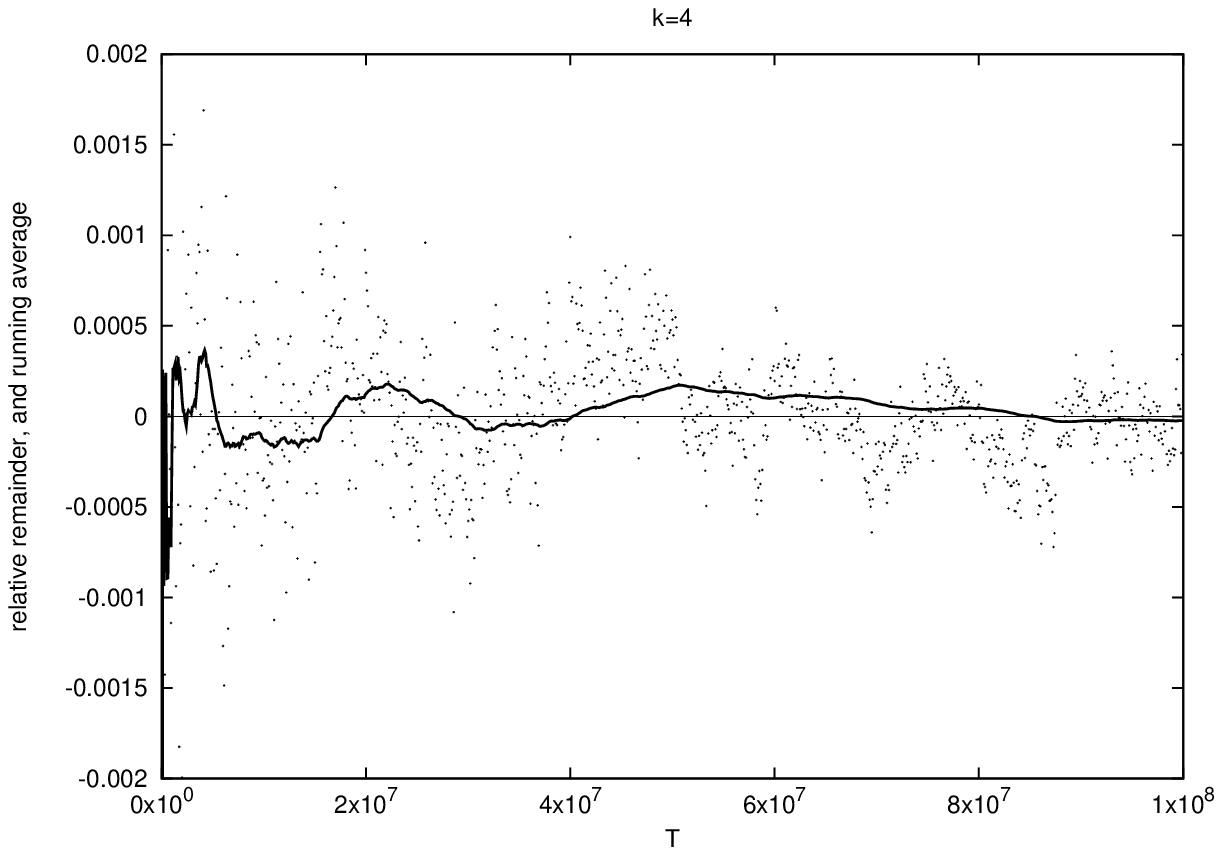}
    }
    \centerline{
       \includegraphics[width=.58\textwidth,height=2.7in]{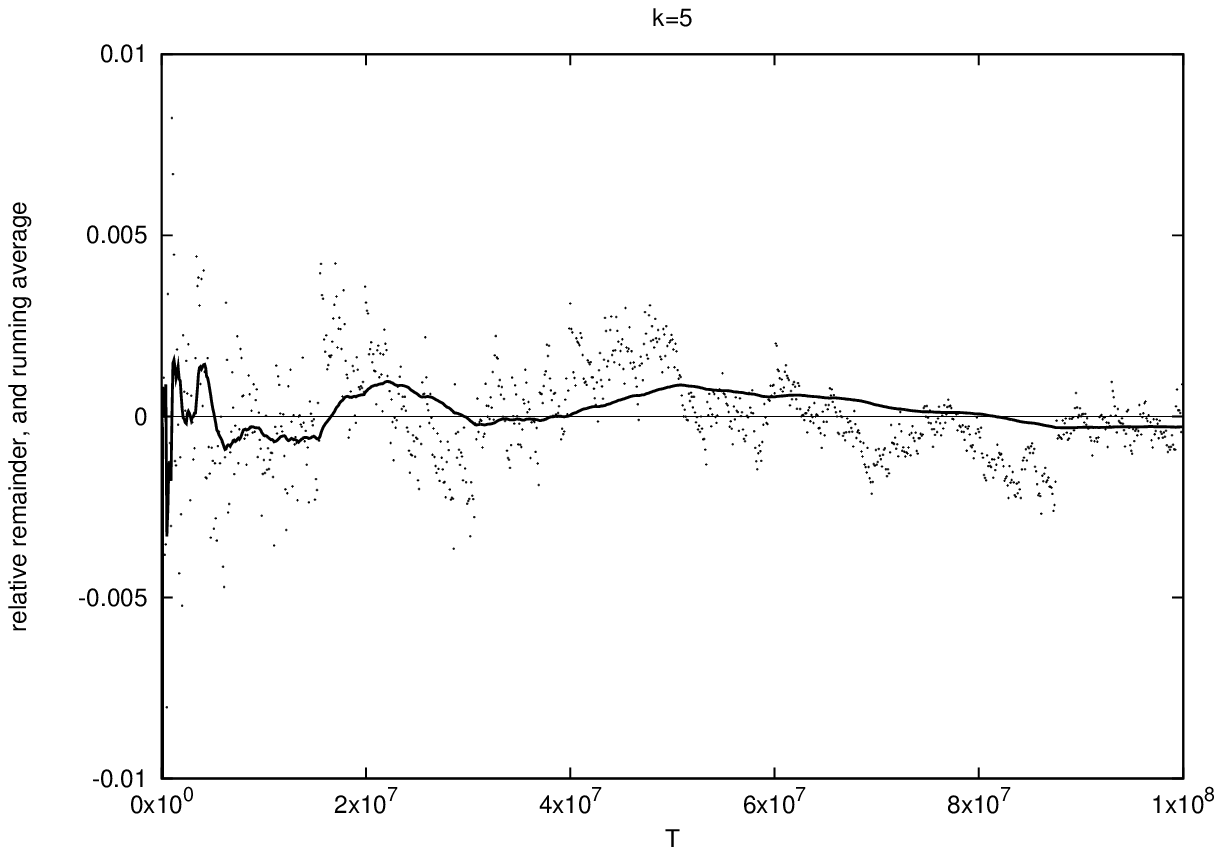}
       \includegraphics[width=.58\textwidth,height=2.7in]{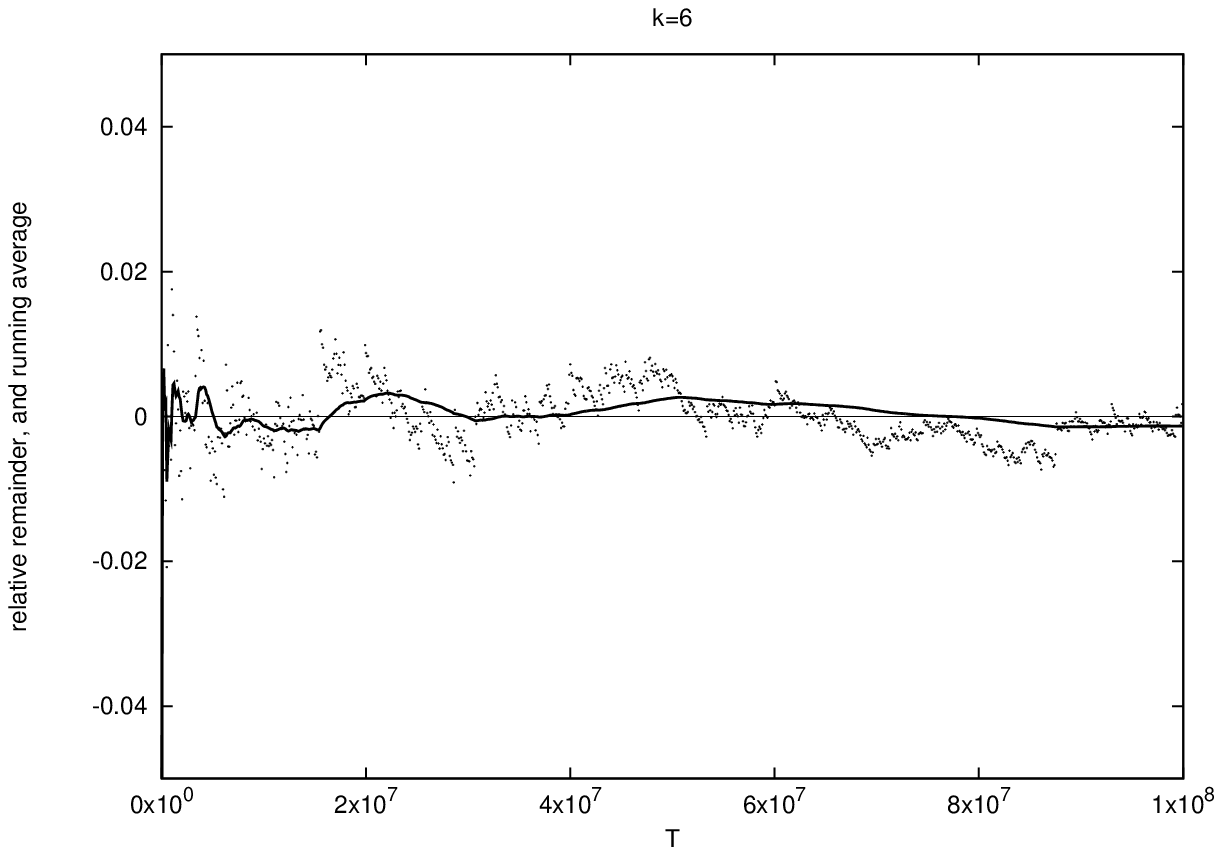}
    }
    \caption{Plot of the relative remainder \eqref{eq:relative remainder}, depicted as dots,
    and running average (solid curve)
    of the remainder \eqref{eq:smooth relative remainder}, for 1000 values of
    $T$, and $1\leq k\leq 6$. Horizontal axis is $T$.}
    \label{fig:1}
\end{figure}

For example, when $k=3$ and $T=10^8$, the $r=2$ term of $P_3$ contributes the
most to~\eqref{eq:conjectured moments}. Furthermore, the first four digits
of~\eqref{eq:conjectured moments}, for $k=3$ and $T=10^8$, arise from the all
the terms with $0 \leq r \leq 8$, while the $r=9$ term contributes to around
the 6th decimal place. Thus, the impact of omitting any of the terms $0\leq r
\leq 8$ would be readily seen in the quality of the comparison in
Table~\ref{table:final_T} or in the plot in Figure~\ref{fig:1}, at least for
$k=3$. For $k=4$ and $t=10^8$, the terms $0 \leq r \leq 12$ all contribute to
the first four digits of~\eqref{eq:conjectured moments}, and, again the quality
of the fit to the actual moment would be noticeably worse if any of these terms
were omitted. The numerics for $k=3$ and $4$ alone provide substantial evidence
favouring the full asymptotics of the moments, being sensitive to almost all
the terms of the conjecture.

Our numerics for larger $k$ also support the conjectured moments, though are
sensitive to a smaller set of terms of $P_k$ for the value $T$ examined. For
example, when $k=7$, the $r=26$ term of $P_k$ contributes the most
to~\eqref{eq:conjectured moments}. Furthermore, the first 2-3 digits
of~\eqref{eq:conjectured moments}, for $k=7$ and $T=10^8$, arise from the terms
with $19 \leq r \leq 33$, whereas the leading coefficient only contributes to
around the 20th decimal place. This illustrates the importance, when testing the
CFKRS prediction, of incorporating the lower terms of $P_k(x)$.

The paper of Hiary and Odlyzko~\cite{HO} contains additional numerics
concerning the moments of zeta. Using values of $c_r(k)$ computed earlier
in~\cite{CFKRS}, they examined the full asymptotics of the moments of
$|\zeta(1/2+it)|$, for $k\leq 6$, but for intervals $T \leq t \leq T+H$. Their
values of $T$ are of size around $10^7$, $10^8$, $10^{15}$, $10^{19}$, and
$10^{22}$, while their values of $H$, are, for the most part, significantly smaller than
$T$ or even small in comparison to $T^{1/2}$.

Their data shows, for intervals with $H$ quite small in comparison to
$T^{1/2}$, a large amount of variation in the computed moments over various
intervals of equal length $H$, whereas the predicted moments change very slowly
as one varies $T$. Presumably one needs $H \gg T^\beta$ for some $\beta>1/2$ in
order to get a good agreement between the computed and predicted moments.
Indeed, their data set `s8', which is comparable to the interval of length
approximately $10^8$ that we examined, shows excellent agreement between
computed moments and predicted moments. Their dataset $z16$, with $T\approx
10^{15}$ and $H\approx 10^8$ shows moderate agreement.

\begin{table}[h!tb]
\centerline{
\begin{tabular}{|c||c|c|c|c|c|c|c|c|c|c|c|c|c|}
\hline
$k$ & 1& 2 & 3 & 4 &5 &6 &7 &8 &9 &10 &11 &12 &13 \cr
\hline
$r$ & 0 & 0 & 2 & 6 & 11 & 18 & 26 & 37 & 49 & 64 & 80 & 99 & 119 \cr
\hline
\end{tabular}
}
\caption{The value of $r$ for each $k$ for which the corresponding term
in~\eqref{eq:conjectured moments} contributes the dominant amount when
$T=10^8$.}
\label{table:max}
\end{table}

To calculate~\eqref{eq:actual moments},
we used the tanh-sinh quadrature scheme~\cite{Ba,BLJ} to accurately
estimate each integral between consecutive non-trivial zeros of the zeta
function on the critical line. All our computations of~\eqref{eq:actual
moments} were carried out using 64 bit machine doubles.
To tabulate all the zeros up to, and slightly beyond, $T=10^8$,
we used Rubinstein's {\tt C++} $L$-function package lcalc~\cite{R}. It applies
the Riemann Siegel formula to evaluate $\zeta(1/2+it)$ and look for sign
changes of the Hardy $Z$-function, Brent's method to compute the zeros of
zeta~\cite{Br}, and a variant of Turing's test~\cite{E} to confirm that
all zeros up to given height have been found.

\begin{figure}[h]
    \centerline{
       \includegraphics[width=.58\textwidth,height=2.7in]{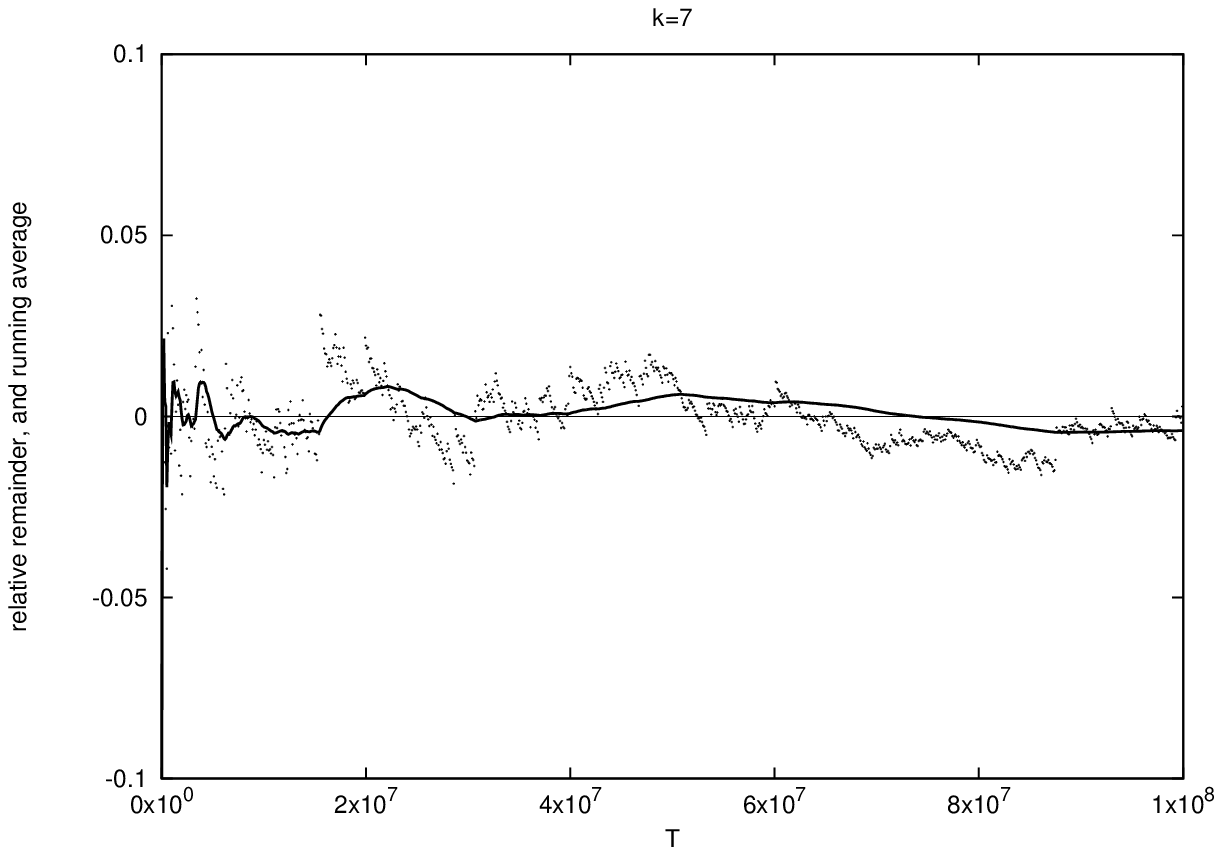}
       \includegraphics[width=.58\textwidth,height=2.7in]{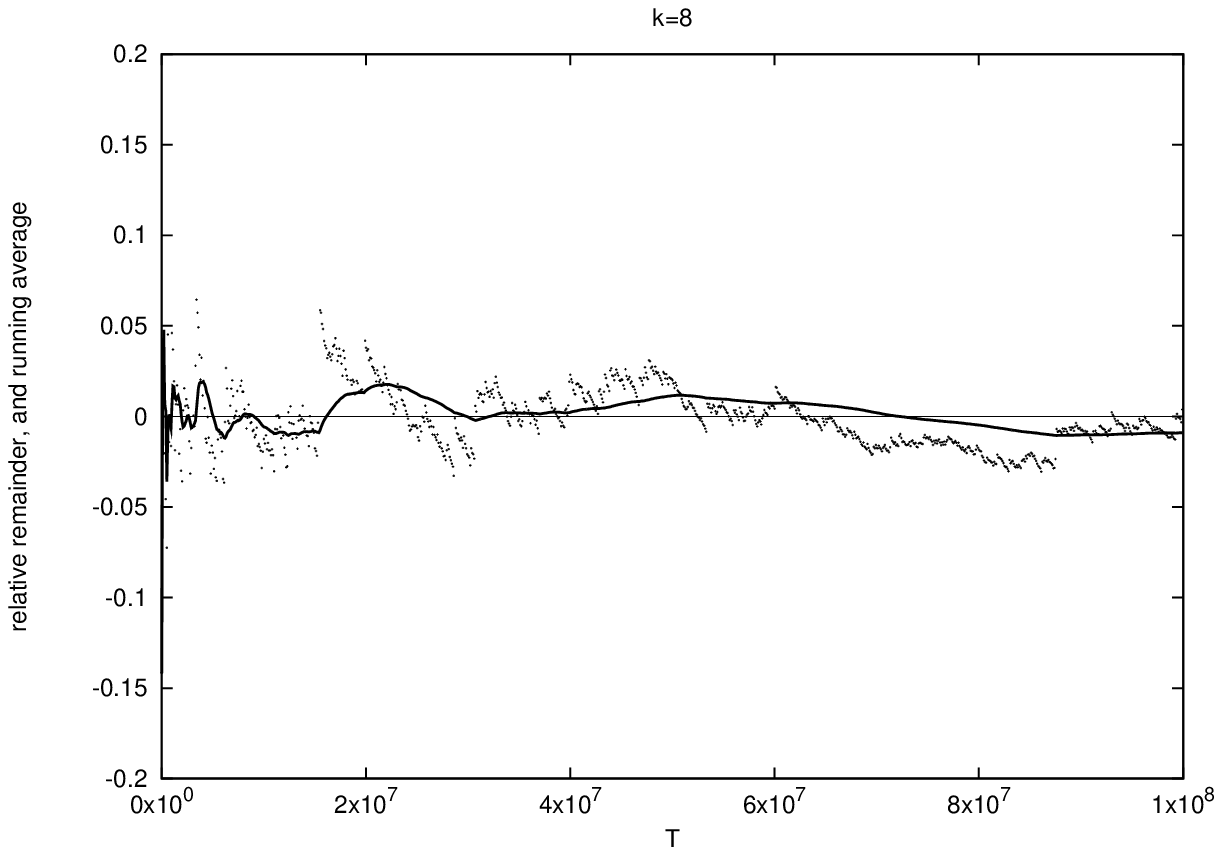}
    }
    \centerline{
       \includegraphics[width=.58\textwidth,height=2.7in]{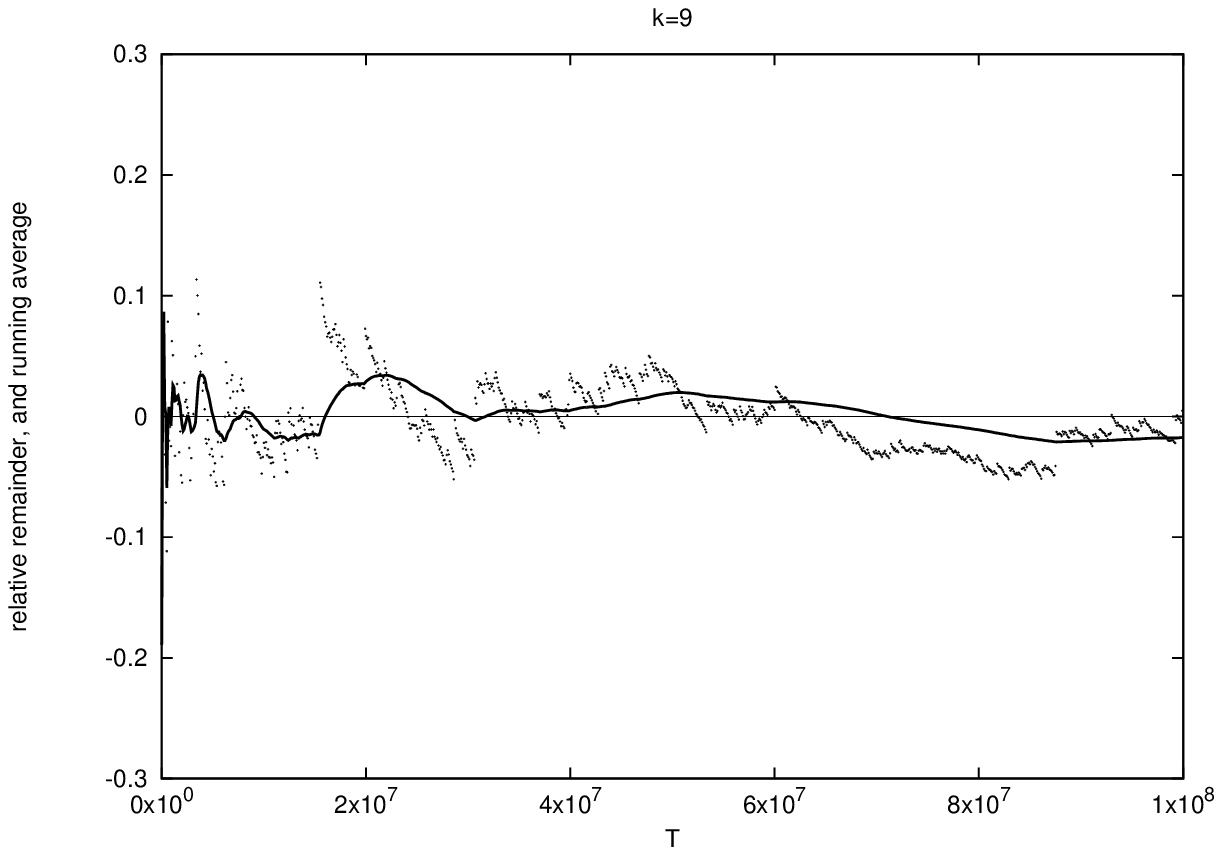}
       \includegraphics[width=.58\textwidth,height=2.7in]{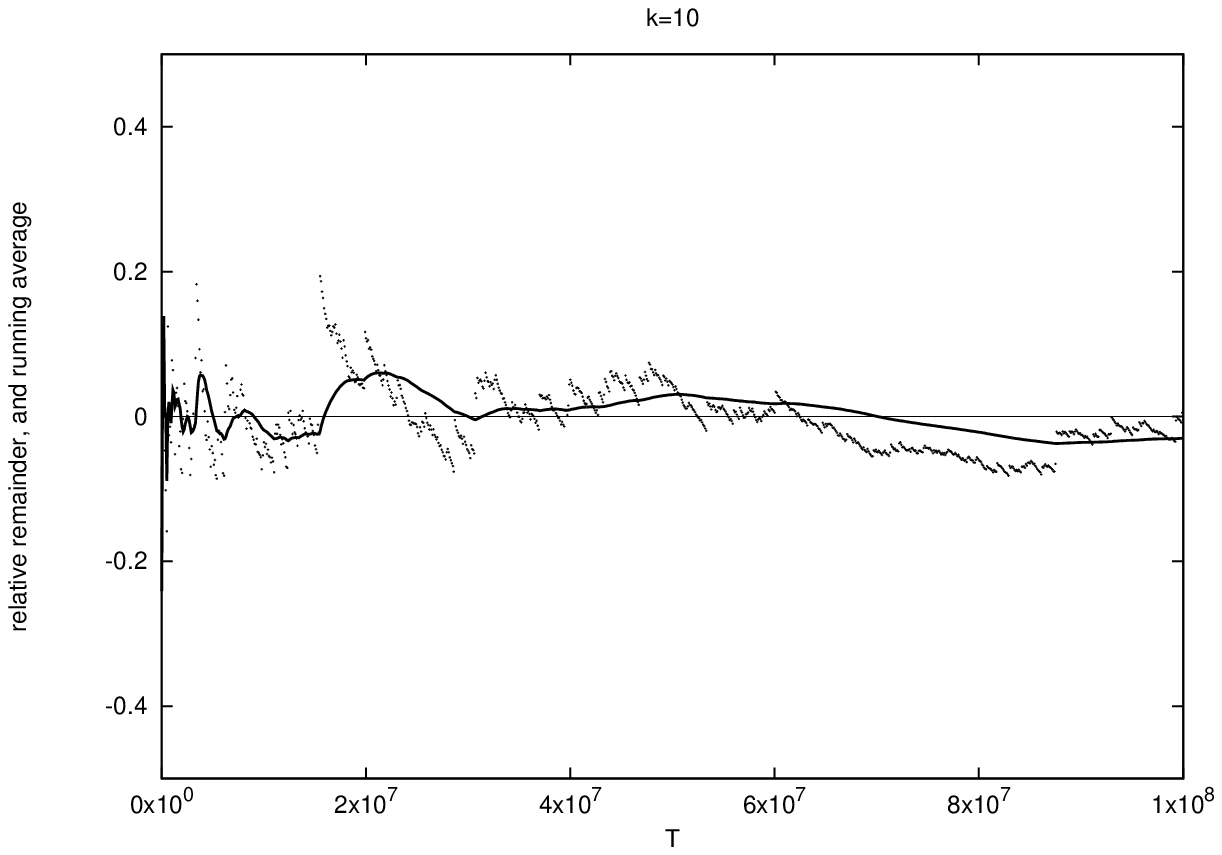}
    }
    \centerline{
       \includegraphics[width=.58\textwidth,height=2.7in]{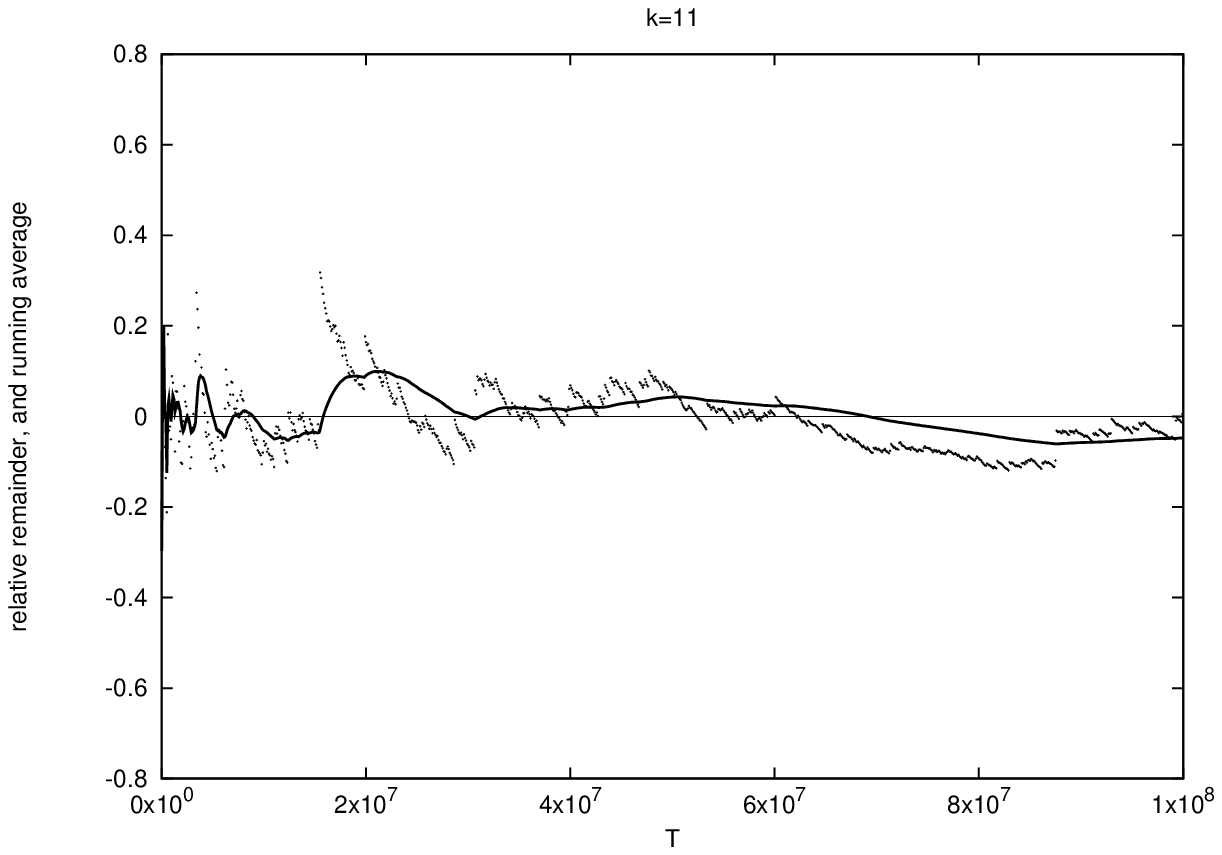}
       \includegraphics[width=.58\textwidth,height=2.7in]{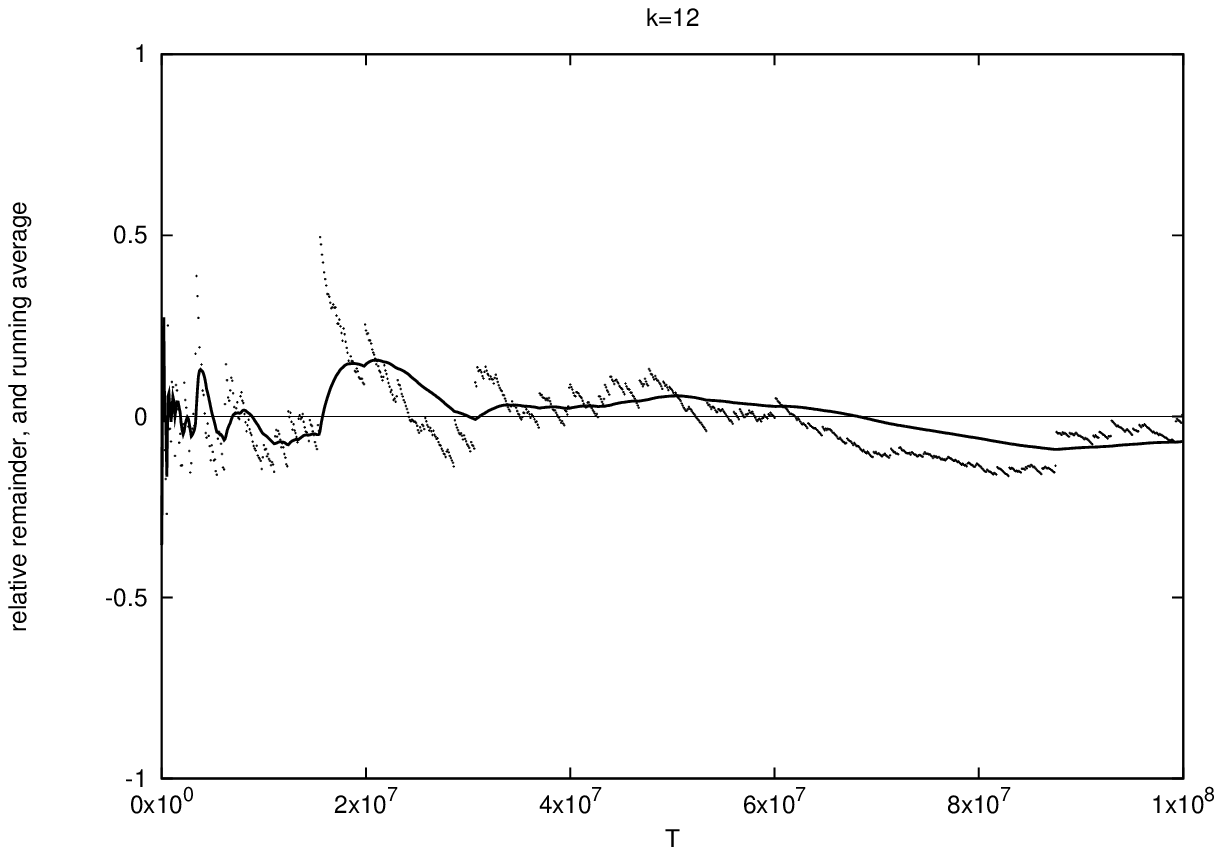}
    }
    \caption{Plot of the relative remainder \eqref{eq:relative remainder}, and running average (solid curve)
    of the remainder \eqref{eq:smooth relative remainder}, for 1000 values of
    $T$, and $7\leq k\leq 12$. Horizontal axis is $T$.}
    \label{fig:2}
\end{figure}

Figures \ref{fig:1}-\ref{fig:2} depict the relative remainder term for 1000 values of
$T$ between 0 and $10^8$, roughly spaced apart every $10^5$.
More specifically, we
let $T_j$ be the imaginary part of the first zero above
$10^5 j$, so that $T_j \approx 10^5j$, and plot the values of
\begin{equation}
    \label{eq:relative remainder}
    \frac{\text{Data}_k(T_j)-\text{Conjecture}_k(T_j)}
    {\text{Conjecture}_k(T_j)},
\end{equation}
for $1 \leq j \leq 1000$, and $1\leq k \leq 12$.

We also calculated the running average of the remainder term divided by the
running average of the conjecture:
\begin{equation}
    \label{eq:running average}
    \frac{\sum_{j=1}^J (\text{Data}_k(T_j)-\text{Conjecture}_k(T_j))}
    {\sum_{j=1}^J \text{Conjecture}_k(T_j)}.
\end{equation}
If we define
\begin{equation}
    \label{eq:smooth actual moments}
    \text{SmoothData}_k(T)= \int_{T_0}^\text{T} |\zeta(1/2+it)|^{2k} (1-t/T) dt,
\end{equation}
and
\begin{equation}
    \label{eq:smooth conjectured moments}
    \text{SmoothConjecture}_k(T)= \int_{T_0}^\text{T} P_k(\log \frac{t}{2 \pi})(1-t/T) dt,
\end{equation}
then~\eqref{eq:running average} gives a discrete approximation to the smoothed
relative remainder:
\begin{equation}
    \label{eq:smooth relative remainder}
    \frac{\text{SmoothData}_k(T_j)-\text{SmoothConjecture}_k(T_j)}{\text{SmoothConjecture}_k(T_j)}.
\end{equation}
The reason for considering smoothed moments is that the noisy remainder terms of the
sharply truncate moments, when averaged, tend to be smaller.

Note that the vertical axes in these figures change from plot to plot to allow us to meaningfully
display the relative remainder terms, which as a whole get larger, as $k$ increases.
We also set the zoom level to
show the running averages in some detail. As a compromise, a few outliers are omitted
from these plots for smaller $T$, roughly up to $10^7$, and $k\leq 4$.

Table~\ref{table:deviation} lists the
standard deviations of 900 values of the remainder term and smoothed remainder term,
for $1\leq k \leq 13$. Specifically, we computed the standard deviation for the values
of~\eqref{eq:relative remainder}, $101\leq j\leq 1000$, and of~\eqref{eq:smooth relative remainder},
for $101 \leq J \leq 1000$.

Table~\ref{table:final_T} lists the values of $\text{Data}_k(T)$, and $\text{Conjecture}_k(T)$
for $k=1,\ldots,13$ and $T=100000000.64$, the first zero of zeta above $10^8$.
We also list the values of the averages over all 1000 values of $T_j$
\begin{equation}
    \label{eq:average data}
    \frac{1}{1000} \sum_{j=1}^{1000} \text{Data}_k(T_j)
\end{equation}
and
\begin{equation}
    \label{eq:average conj}
    \frac{1}{1000} \sum_{j=1}^{1000} \text{Conjecture}_k(T_j).
\end{equation}

Our data supports the CFKRS conjecture for the full asymptotics of
the moments of zeta as  described in equation~\eqref{eq:moments zeta}, though
for larger $k$, it is difficult to gauge the size of the remainder term.

For $k=1$ the data suggests an even stronger remainder term of
$O\left(T^{1/4+\delta}\right)$, supported by the agreement of $\text{Data}_1(T)$ with
$\text{Conjecture}_1(T)$ to about $3/4$ of the decimal places left of the
decimal point for the values of $T$ examined. The relative remainder term is
of size around $10^{-6}$ when $T \approx 10^8$, and, typically,
an order of magnitude smaller when averaged.

For $k=2$ the data agrees with the conjecture to about half the decimal places,
with a relative remainder term of size around $10^{-6}$. For $k=3$, the
agreement is to slightly less than half the decimal places.

For fixed $T$, as $k$ increases, the moments have the effect of amplifying the
largest values of $|\zeta(1/2+it)|$. This can be seen in our plots, for
larger $k$, where the remainder terms are qualitatively the same, with large jumps at
the same values of $T$ corresponding to relatively large values of zeta. It
therefore becomes more difficult to ascertain, as $k$ grows, whether an upper
bound for the remainder term of the form $O_k(T^{1/2+\delta})$ holds.
Nonetheless, the table and figures reveal an excellent fit for the CFKRS
prediction with the moments which persists through to the $24$th and $26$th
moments, where the relative agreement is to within around one to two decimal places.

In some sense, the fit between columns 2 and 3 of Table~\ref{table:final_T} is
better than it ought to be for larger $k$, for example
more than three decimal places for $k=13$, but only agreeing to around 90\%
for the running average of the remainder.
A quick inspection of the figures reveals that the relative remainder, at $T=10^8$,
happens to, fortuitously, best its neighbours, especially for larger $k$.
Nonetheless, the overall agreement
between the CFKRS prediction and our data across all values of $k$ and $T$,
as depicted in the figures, lends strong support to their conjecture for
the full asymptotics of the moments of the zeta function.

\begin{table}[h!tb]
\centerline{
\begin{tabular}{|c|c|c||c|c|}
\hline
$k$ & $\text{Data}_k(T)$ & $\text{Conjecture}_k(T)$ & \eqref{eq:average data} & \eqref{eq:average conj} \cr
\hline
1 & 1673723690.436 & 1673723498.495 & 62463107.03367 & 62463106.44834 \cr
2 & 637388343407. & 637389923500. & 22091815715.8 & 22091815007.3 \cr
3 & 8.04585314342e+14 & 8.04581403344e+14 & 2.54969941363e+13 & 2.54969410053e+13 \cr
4 & 1.7376480696e+18 & 1.7374512576e+18 & 5.0233293703e+16 & 5.0234406051e+16 \cr
5 & 5.0837678819e+21 & 5.0816645028e+21 & 1.3436072658e+20 & 1.3439894636e+20 \cr
6 & 1.815301994e+25 & 1.813639687e+25 & 4.400401104e+23 & 4.406152513e+23 \cr
7 & 7.480512969e+28 & 7.468884126e+28 & 1.66796929e+27 & 1.674460675e+27 \cr
8 & 3.43851173e+32 & 3.43090327e+32 & 7.0665341e+30 & 7.13020297e+30 \cr
9 & 1.72388578e+36 & 1.71918466e+36 & 3.26838252e+34 & 3.32617445e+34 \cr
10 & 9.2785049e+39 & 9.251733e+39 & 1.6228859e+38 & 1.6729912e+38 \cr
11 & 5.2991086e+43 & 5.2863072e+43 & 8.5447722e+41 & 8.967058e+41 \cr
12 & 3.182548e+47 & 3.179455e+47 & 4.726347e+45 & 5.076008e+45 \cr
13 & 1.995625e+51 & 1.999377e+51 & 2.726982e+49 & 3.013406e+49 \cr
\hline
\end{tabular}
}
\caption{The values of $\text{Data}_k(T)$ and $\text{Conjecture}_k(T)$ at $T=100000000.64$,
the first zero after $10^8$. We also display, in the last two columns, the values of
the averages over all 1000 values of $T_j$:
$\frac{1}{1000} \sum_{j=1}^{1000} \text{Data}_k(T_j)$, and
$\frac{1}{1000} \sum_{j=1}^{1000} \text{Conjecture}_k(T_j)$.
}
\label{table:final_T}
\end{table}

\begin{table}[h!tb]
\centerline{
\begin{tabular}{|c|c|c|}
\hline
$k$ & $\text{Standard deviation}$ & $\text{Standard deviation, smoothed}$ \cr
\hline
1 & 5.0337e-07 & 4.4433e-08 \cr
2 & 7.6993e-06 & 2.7833e-07 \cr
3 & 5.9228e-05 & 6.5080e-06 \cr
4 & 0.00032311 & 8.9497e-05 \cr
5 & 0.0012292 & 0.00046265 \cr
6 & 0.0034752 & 0.0015118 \cr
7 & 0.0079677 & 0.0038227 \cr
8 & 0.015741 & 0.0081774 \cr
9 & 0.027867 & 0.015496 \cr
10 & 0.045375 & 0.026757 \cr
11 & 0.069216 & 0.042914 \cr
12 & 0.10024 & 0.064832 \cr
13 & 0.1392 & 0.093254 \cr
\hline
\end{tabular}
}
\caption{The standard deviations of 900 values of the remainder terms, and smoothed remainder terms,
for $1\leq k \leq 13$. Specifically, we computed the standard deviation for the values
of~\eqref{eq:relative remainder}, $101\leq j\leq 1000$, and of~\eqref{eq:smooth relative remainder},
for $101 \leq J \leq 1000$.
}
\label{table:deviation}
\end{table}

{\bf Acknowledgements:} We would like to thank the referee for his insightful feedback.

\bibliographystyle{amsplain}

\begin{thebibliography}{10}


\bibitem[Ba]{Ba}
David H. Bailey, {\it Tanh-Sinh High-Precision Quadrature}, Jan 2006; LBNL-60519.

\bibitem[BLJ]{BLJ}
David H. Bailey, Xiaoye S. Li and Karthik Jeyabalan, {\it A Comparison of Three High-Precision Quadrature Schemes}, Experimental Mathematics, vol. 14 (2005), no. 3, pg 317-329. LBNL-53652.

\bibitem[Br]{Br} Richard P. Brent, {\it Algorithms for Minimization without Derivatives},
Chapter 4, Prentice-Hall, Englewood Cliffs, NJ, 1973.

\bibitem[CFKRS]{CFKRS}
J. B. Conrey, D. W. Farmer, J. P. Keating, M. O. Rubinstein, and N. C. Snaith,
{\it  Integral moments of ${L}$-functions},
\newblock Proceedings of the London Mathematical Society (3), 91 (2005), 33--104.

\bibitem[CFKRS2]{CFKRS2}
J. B. Conrey, D. W. Farmer, J. P. Keating, M. O. Rubinstein, and N. C. Snaith,
{\it Lower order terms in the full moment conjecture for the Riemann Zeta Function},
Journal of Number Theory, Volume 128, Issue 6, June 2008, 1516-1554.

\bibitem[DGH]{DGH}
\newblock   A. Diaconu, D. Goldfeld, and J. Hoffstein,
Multiple Dirichlet series and moments of zeta- and $L$-functions,
Compositio Math. 139 (2003), no.~3, 297-360.

\bibitem[E]{E}
\newblock 
H. Edwards, {\it Riemann's Zeta Function}, Academic Press, 1974.

\bibitem[FHLPZ]{FHLPZ}
\newblock
L. Fousse, G. Hanrot, V. Lefevre, P. P\'elissier, P. Zimmermann,
{\it MPFR: A multiple-precision binary floating-point library with correct rounding},
ACM Transactions on Mathematical Software, (2), 33 (2007).

\bibitem[H-B]{H-B}
\newblock
D. R. Heath-Brown,  {\it The fourth power moment of the Riemann zeta-function},
Proc. London Math. Soc. (3) 1979 {\bf 38}  pp. 385 -- 422.

\bibitem[HO]{HO}
G. A. Hiary and A. M. Odlyzko, {The zeta function on the critical line:
numerical evidence for moments and random matrix theory models}, Math.
Comp., 81 (2012), 1723--1752.

\bibitem[HR]{HR}
G. A. Hiary and  M. O. Rubinstein, {\it Uniform asymptotics for the full moment
conjecture of the Riemann zeta function}, Journal of Number Theory, Volume 132,
Issue 4, (2012), 820--868.

\bibitem[I]{I}
\newblock
 A. E. Ingham,  {\it Mean-value theorems in the theory of the Riemann
 zeta-function},   Proceedings of the London Mathematical Society  (92) 1926
 {\bf 27}  pp. 273--300.

\bibitem[R]{R} M. O. Rubinstein, {\it lcalc - the $L$-function calculator},
{\tt http://code.google.com/p/l-calc/}.

\end{thebibliography}

\end{document}